\documentclass{elsart}  
 
\usepackage{fullpage} 
 
\usepackage{amssymb} 
\usepackage{amsmath} 
\usepackage{amsfonts} %
\newcommand{\baseRing}[1]{\ensuremath{\mathbb{#1}}}

\newcommand{\R}{\baseRing{R}} 
 
\newcommand{\C}{\baseRing{C}}  
 
\newcommand{\trace}{\mathop{\mathrm{trace}}} 
\newcommand{\Aut}{\mathop{\mathrm{Aut}}} 
\newtheorem{theorem}{Theorem}[section] 
\newtheorem{lemma}[theorem]{Lemma}

\newtheorem{definition}[theorem]{Definition} 
\newtheorem{remark}[theorem]{Remark} 
\newtheorem{example}[theorem]{Example}

\begin{document} 
\begin{frontmatter} 
\title{Symmetry groups, semidefinite programs, \\ and sums of squares} 
\author[label1]{Karin Gatermann} 
\author[label2]{Pablo A. Parrilo} 
\address[label1]{Konrad-Zuse-Zentrum f\"ur 
Informationstechnik, Takustr.\ 7, D-14195 Berlin, Germany or Institut~I,  
Fachbereich Mathematik und Informatik, FU Berlin, Arnimallee 2-6, 
D-14195 Berlin, Germany} 
\ead{gatermann@zib.de} 
\ead[url]{http://www.zib.de/gatermann} 
\address[label2]{Automatic Control Laboratory,  
Physikstrasse 3 / ETL, ETH-Zentrum, CH-8092 Z\"urich, Switzerland} 
\ead{parrilo@aut.ee.ethz.ch} 
\ead[url]{http://www.aut.ee.ethz.ch/\~{}parrilo} 
\begin{abstract} 
  We investigate the representation of symmetric polynomials as a sum
  of squares. Since this task is solved using semidefinite programming
  tools we explore the geometric, algebraic, and computational
  implications of the presence of discrete symmetries in semidefinite
  programs.  It is shown that symmetry exploitation allows a
  significant reduction in both matrix size and number of decision
  variables.  This result is applied to semidefinite programs arising
  from the computation of sum of squares decompositions for
  multivariate polynomials. The results, reinterpreted from an
  invariant-theoretic viewpoint, provide a novel representation of a
  class of nonnegative symmetric polynomials.  The main theorem states
  that an invariant sum of squares polynomial is a sum of inner
  products of pairs of matrices, whose entries are invariant
  polynomials. In these pairs, one of the matrices is computed based
  on the real irreducible representations of the group, and the other
  is a sum of squares matrix. The reduction techniques enable the
  numerical solution of large-scale instances, otherwise
  computationally infeasible to solve.

\noindent 
{\bf AMS Subject Classification:} 14P99, 26C05, 13A50, 68W30, 90C22 
\end{abstract} 
\begin{keyword} real algebraic geometry, convex optimization, finite groups,  
symmetry 
\end{keyword} 
\end{frontmatter}

\section{Introduction and preliminaries} 
\label{sec:intro} 
 
A fundamental problem in real algebraic geometry is the existence and 
computation of a representation of a multivariate polynomial as a sum 
of squares (SOS). In other words, the question of finding $p_i\in 
\R[{\bf x}],i=1,\ldots,N$ such that 
\[ 
f(\mathbf{x})=\sum_{i=1}^N \bigl(p_i(\mathbf{x})\bigr)^2. 
\] 
Besides its obvious theoretical interest, this is also a relevant 
question in many areas of applied mathematics, such as continuous and 
combinatorial optimization, as well as real polynomial system solving 
via the connections with the Positivstellensatz 
\cite{pablo,sdprelax,Sturm02}. While algebraic techniques have been 
proposed to decide the problem \cite{PowWor}, recent results 
\cite{shor,pablo,PaStu01} suggest that this problem can be solved much 
more efficiently using numerical optimization techniques such as 
semidefinite programming (SDP) \cite{HandSDP}. 
 
It is a fact that many problems in applied mathematics and engineering 
possess inherent symmetries in their mathematical formulations. While 
in some cases this mirrors the underlying structure of existing 
physical systems, these features can also arise as a result of the 
chosen mathematical abstraction. A natural question, therefore, is the 
investigation of the consequences of these symmetries in the 
efficiency and performance of the algorithms associated with the 
computation of diverse quantities of interest.  
 
Motivated by these reasons, we study the problem of finding sum of 
squares decompositions whenever the given polynomial $f$ is invariant 
under the action of a \emph{finite group} (discrete symmetries).  Our 
initial focus will be on the consequences of symmetry in general 
semidefinite programs, emphasizing later the specific class that 
appears in the computation of sum of squares decompositions and the 
associated Positivstellensatz relaxations. 
 
In general, symmetry reduction techniques have been explored in 
several different contexts.  A partial list of relevant examples and 
application areas are dynamical systems and bifurcation theory 
\cite{Golubitsky}, polynomial system solving \cite{Wor94,gat}, 
numerical treatment of partial differential equations \cite{Fassler}, 
and Lie symmetry analysis in geometric mechanics \cite{MechSym}. 
Despite belonging to very different domains, common themes in all 
these developments are the conceptual and computational benefits of 
a thorough understanding of the constraints imposed by symmetry.   
 
The general advantages of symmetry exploitation are numerous: 
 
\begin{description} 
\item [Problem size.] The first immediate advantage is the reduction 
  in problem size, as the new instance can have a significantly 
  smaller number of variables and constraints. This notably affects 
  the total running time. A well-known application area where these 
  techniques are employed is finite-element methods for partial 
  differential equations \cite{Fassler}. 
   
\item [Degeneracy removal.] A very common source of difficulties in 
  numerical methods is the presence of repeated eigenvalues of high 
  multiplicity. These usually arise from structural properties, and 
  can be removed by a proper handling of the symmetry. 
   
\item [Conditioning.] For the above mentioned reasons, the resulting 
  symmetry-aware methodologies in general have much better numerical 
  conditioning. Additionally, the elimination of non-essential 
  variables allows for faster and more accurate solutions. 
   
\item [Reliability.] Smaller size instances are usually less prone to 
  numerical errors. For instance, numerical path following is much 
  more accurate in a fixed point space of smaller dimension than in 
  the original high dimensional space. 
\end{description} 
In later sections, we illustrate the concrete form that these general 
benefits take in the problem of computation of sum of squares 
decompositions. 
 
For the exploitation of symmetry of finite groups there are two 
different, but complementary, approaches. The first one is the theory 
of \emph{linear representations}, as presented in 
\cite{Fassler,Serre}, and the second one being \emph{invariant theory} 
\cite{sturmfelsinvariant,DerKem02}.  We make use of both approaches in 
our work, and show how the results can be naturally interpreted in 
either framework.

A particularly interesting application of the sum of squares (SOS) 
methods, outlined in Section~\ref{sec:sossdp}, is the optimization of 
symmetric polynomials.  In \cite{pablo,PaStu01}, and based on 
\cite{shor}, an SDP-based algorithm for the computation of a global 
lower bound for polynomial functions is described. When the polynomial 
to be minimized has certain symmetry properties, then the theory 
described here can be fruitfully applied.  In related work, Choi 
\emph{et al.} in \cite{lam} have shown how symmetry reduction 
techniques can significantly reduce the complexity of analyzing 
SDP/SOS conditions. 
 
The outline of the paper is as follows: in Section~\ref{sec:sossdp} we
start by recalling the intimate connections between sum of squares and
semidefinite programs, as well as their applications in optimization.
In Section~\ref{sec:sdp} the concept of an \emph{invariant
semidefinite program} is introduced.  The key result in this section
relies critically on the convexity property, and shows that the
feasible set can always be restricted to a specific invariant
subspace, thus reducing the problem to a semidefinite program of
smaller dimensions.  In Section~\ref{sec:blk}, linear representation
theory is employed to find a change of coordinates where all the
matrices in the SDP have block diagonal form. The SDP therefore
collapses into a collection of coupled smaller optimization problems,
which are much easier to solve.  In the following section, we
concentrate on the sum of squares techniques for symmetric
polynomials, where the group representation on matrices is induced by
an action on the space of monomials. We show that the complexity of
checking positive definiteness conditions simplifies dramatically, and
illustrate the techniques with numerous examples.  Examples from the
literature \cite{choiform,Reznick,PaStu01,pablo,Sottile} have been
computed with SOSTOOLS \cite{SOStools}.  An invariant-theoretic
viewpoint is introduced in Section~\ref{sec:invar}, stressing its
computational efficiency and the natural connections with the
representation-based approach developed earlier. As part of the
developments, we obtain a novel and very appealing representation of
symmetric SOS polynomials as a sum of inner products between invariant
matrices, one of them being an SOS matrix (Theorem
\ref{thm:matrixsosrep}). Finally, in Section~\ref{sec:savings}, the
computational savings of our techniques are illustrated in some
large-scale problems.

\section{Optimization of symmetric polynomials  
and sum of squares decompositions} 
\label{sec:sossdp} 

It is well-known that the problem of global minimization for 
multivariate polynomials is an NP-hard problem. The technique we 
describe, which was first proposed by Shor \cite{shor}, is a 
relaxation: it produces a certifiable lower bound for the optimal 
value, and in practice very often provides the correct answer. We 
provide a brief outline below. 
 
Consider a given polynomial $f(x_1,\ldots,x_n)$ of even total degree 
$2d$. If we can find a real number $\lambda$ and polynomials $p_i 
\in \R[x_1,\ldots,x_n], i=1,\ldots,N$ such that: 
\begin{equation} 
f(\mathbf{x}) - \lambda = \sum_{i=1}^N p_i(\mathbf{x})^2,  
\label{eq:sosdecomp} 
\end{equation} 
then clearly $\lambda$ is a global lower bound of $f(\mathbf{x})$, 
i.e., $f(\mathbf{x}) \geq \lambda$ for all $\mathbf{x} \in \R^n$. As 
shown in \cite{shor}, such a decomposition, if it exists, can be found 
using convex optimization, and in particular, semidefinite programming 
(see \cite{PaStu01,N98polynomials,Lasserre} and the references therein). 

We illustrate the idea in the general case, i.e., when the polynomial 
$f$ is dense. We write $f(\mathbf{x})=Y^TAY$, where $A$ is a real 
symmetric matrix and $Y$ denotes the column vector whose entries are 
all the monomials in $\,x_1,\ldots,x_n$ of degree at most $d$.  The 
length of the vector $Y$ is equal to the binomial coefficient 
\[ 
 N \,\, = \,\, \binom{n+d}{d} . 
\] 
However, the matrix $A$ in this quadratic representation is not 
necessarily unique.  Let $\mathcal{S}^N$ denote the vector space of 
real symmetric $N\times N$-matrices. Then 
\[ 
\mathcal{L}_f=\{A\in \mathcal{S}^N\,\,|\,\, f(\mathbf{x})=Y^TAY\}. 
\] 
denotes the set of all real symmetric $N \times N$-matrices $A$ for 
which the representation is valid.  This is an affine subspace in the 
space of real symmetric $N \times N$-matrices.  Assume that the 
constant monomial $1$ is the first entry of $Y$.  Let $E_{11}$ denote 
the matrix unit whose only non-zero entry is a one in the upper left 
corner.  Then $\mathcal{L}_{f-\lambda}= \{X\in \mathcal{S}^N\,\,|\,\, 
f(\mathbf{x})-\lambda=Y^TQY\}$ where $Q = A -\lambda \cdot E_{11}$. 
 
\begin{lemma} (\cite{PaStu01}) 
\label{lemma1} Let $f(\mathbf{x})\in \R[{\bf x}]$.  
For any real number $\lambda$, the following statements are equivalent: 
\begin{enumerate} 
\item[(i)] The polynomial $ \, f({\bf x}) - \lambda \,$ has a 
  representation as a sum of squares in $\R[{\bf x}]$. 
\item[(ii)]  
There is a matrix $ Q \in \mathcal{ L}_{f-\lambda}$ such that 
$Q$ is positive semidefinite, 
that is, all eigenvalues of this symmetric matrix 
are nonnegative reals. 
\item[(iii)]  
There is  a matrix $ A \in \mathcal{ L}_f$ such that 
$\,A - \lambda \cdot E_{11} \,$ is positive semidefinite. 
\end{enumerate} 
\label{lem:polysdp} 
\end{lemma} 
{\bf Proof:} 
The matrix $\,A - \lambda \cdot E_{11} \,$ is positive semidefinite 
if and only if there exists a real factorization 
$\,  A - \lambda \cdot E_{11} \, = \, B^T \cdot B $.  
If this holds then 
\[ 
\begin{array}{rcl} 
 f({\bf x}) - \lambda & = & 
  Y^T \cdot A \cdot Y - \lambda  \, = \,  
  Y^T \cdot (  A - \lambda \cdot E_{11}) \cdot Y \\ 
& = & 
  Y^T \cdot   B^T B \cdot Y \, = \,  
(B Y)^T \cdot (B Y)=\sum_{i=1}^N\bigl((BY)_i\bigr)^2 
\end{array} 
\] 
is a sum of squares, and every representation as sum of squares 
arises in this way. 
\hfill \qed 
 
As we shall see next, Lemma~\ref{lem:polysdp} implies that we can find 
representations of a polynomial as a sum of squares by solving a 
particular kind of convex optimization problems called 
\emph{semidefinite programs}. 
 
The standard semidefinite program \cite{HandSDP} is given in the 
following way.  The inner product on $\mathcal{S}^N$ is denoted by 
$\langle \cdot , \cdot \rangle$ and let $\mathcal{S}^N_+ \subset 
\mathcal{S}^N$ denote the set of positive semidefinite matrices.  It 
is well-known that the set $\mathcal{S}^N_+$ forms a {\em convex 
  cone}.  Given a symmetric matrix $C\in\mathcal{S}^n $ we define a 
linear cost function $F:\mathcal{S}^N_+\rightarrow \R$ by 
$F(X)=\langle C,X\rangle$.  Given an affine subspace $\mathcal{L}$ of 
$\mathcal{S}^N$, the convex optimization problem is: 
\[ 
 F:=\min_{X \in \mathcal{L} \cap \mathcal{S}^N_+} \langle C , X 
\rangle, 
\] 
where we minimize the linear cost function $F(X)=\langle C,X\rangle$, over  
the set of feasible matrices which is the intersection of an 
affine matrix subspace $\mathcal{L}$ with the set $\mathcal{S}^N_+$ 
of positive semidefinite matrices. 
This is a finite dimensional convex optimization problem, for which 
efficient algorithms are available \cite{HandSDP}. 
 
For the case of sum of squares the linear subspace of symmetric 
matrices representing $f(x)-\lambda$ is $\mathcal{L}=\{X\in 
\mathcal{S}^N\,\,|\,\,\exists\, \lambda\in \R \mbox{ with } 
f(x)-\lambda=Y^TXY\}$ while the cost function is $F(X)=\langle 
C,X\rangle$ with $C=E_{11}$.  We write $f^{sos}$ for the largest real 
number $\lambda$ for which the minimum of $F(X)$ is attained. 
$f^{sos}$ is a lower bound for the global minimum $f^*$ of $f$ on 
$\R^n$. 
This inequality may be strict. It is even possible that  
$\,f^{sos} = - \infty $ for the case that no $\lambda$ exists such that 
$f(x)-\lambda$ has a representation as sum of squares.  We refer the 
reader to \cite{PaStu01} for more background, numerical results, and 
comparison with other methodologies.  
 
{\bf Sum of squares matrices:} We introduce now the notion of
\emph{sum of squares matrices}. These will be matrices, with
polynomial entries, that are ``obviously'' positive semidefinite for
all values of the indeterminates. Therefore, this provides a natural
generalization of the sum of squares condition as an easily verifiable
certificate for polynomial nonnegativity.
\begin{definition} 
Let $S \in \R[\mathbf{x}]^{m \times m}$ be a symmetric matrix, and $\mathbf{y} = [y_1, 
\ldots, y_m]$ be new indeterminates. The matrix $S$ is a \emph{sum of 
squares (SOS) matrix} if the scalar polynomial $\mathbf{y}^T S \mathbf{y}$ 
is a sum of squares in $\R[\mathbf{x},\mathbf{y}]$. 
\label{def:sosmatrix} 
\end{definition} 
It is clear that this concept reduces to the standard SOS notion for 
$1 \times 1$ matrices, i.e., scalar polynomials.  
\begin{example} 
Consider the matrix $S \in \R[x]$, given by: 
\[ 
S = \left[ 
\begin{array}{cc} 
x^2-2x+2 & x \\ x & x^2  
\end{array}\right]. 
\] 
This is a SOS matrix, since 
\[ 
y^T S y = 
\left[\begin{array}{c}
y_1 \\
x y_1 \\
y_2 \\
x y_2 
\end{array}\right]^T
\left[\begin{array}{rrrr}
2 & -1 & 0 & 1 \\
-1 &  1 & 0  &  0 \\
0  & 0   & 0  & 0  \\
1  & 0   & 0  & 1 \\
\end{array}\right]
\left[\begin{array}{c}
y_1 \\
x y_1 \\
y_2 \\
x y_2 
\end{array}\right] =
(y_1 + x \, y_2)^2 + (x \, y_1 -y_1)^2. 
\]\hfill$\diamondsuit$ 
\end{example} 
An important result about SOS matrices is the fact that $S \in \R[x]^{m
\times m}$ (notice $n=1$ here) is positive semidefinite for all $x$ if
and only if it is a SOS matrix. This is a direct restatement of the
fact that $(2,n;m,2)$ positive semidefinite biforms are SOS; see
\cite{CLRrealzeros} and the references therein.

Since the SOS matrix condition is equivalent to a scalar SOS test in a
polynomial with more variables, it is obvious that this can also be
computed using semidefinite programming techniques. Furthermore,
because of the special bipartite structure, only monomials of the form
$x_i^k y_j$ will appear in the vector $Y$, allowing for a more
efficient implementation.

As mentioned, in many practical applications the polynomial $f$ has 
symmetry properties, such as being invariant under a group of 
transformations. 
\begin{example} 
Consider the quartic trivariate polynomial analyzed in 
  \cite{PaStu01}: 
\[ 
f(x,y,z) = x^4 + y^4 + z^4 - 4 \, x y z +x + y + z. 
\] 
This polynomial is invariant under all permutations of $\{x,y,z\}$ 
(the full symmetric group $S_3$).  The optimal value is $f^* = 
-2.112913882$, achieved at the orbit of global minimizers: 
$$            (0.988, -1.102, -1.102) \, , \,\, 
            (-1.102, 0.988, -1.102) \, , \,\, 
            (-1.102, -1.102, 0.988). 
$$ 
For this polynomial, it holds that $f^{sos} = f^*$. 
\hfill$\diamondsuit$ 
\label{ex:berndpoly}
\end{example} 
The symmetry appearing in polynomials (such as $f$ above) suggest that 
it may be possible to \emph{customize} the general SDP-based 
methodology to obtain more efficient algorithms. We present next our 
results for a new class of symmetric semidefinite programs, 
specializing them later to those that arise in the computation of sum 
of squares decompositions of symmetric polynomials. 

\section{Invariant semidefinite programs} 
\label{sec:sdp} 

In order to introduce our techniques for the computation of sum of 
squares decompositions for {\em invariant} polynomials we start by 
analyzing the consequences of symmetry in arbitrary semidefinite 
programs.  We consider first the most general cases, restricting in 
later sections to the special symmetries that arise in our specific 
class of problems.  The major result is that the symmetric 
semidefinite optimization problem can always be reduced to an easier, 
smaller problem because of its {\em convexity}, see 
Theorem~\ref{thm:restricted} below. 
 
As mentioned in the previous section, a semidefinite program is 
defined as follows: given an affine subspace $\mathcal{L}$ of 
$\mathcal{S}^N$ and a symmetric matrix $C\in \mathcal{S}^N$ we 
have 
\begin{equation} 
F:=\min_{X \in \mathcal{L} \cap \mathcal{S}^N_+} \langle C , X 
\rangle, 
\label{eq:sdp1} 
\end{equation} 
where we minimize the linear cost function $F(X)=\langle C,X\rangle$, 
over the set of feasible matrices $\mathcal{L}\cap\mathcal{S}_+^N$. 
Again $\mathcal{S}^N \subset \R^{N \times N}$ denotes the space of 
symmetric $N\times N$ real matrices, equipped with the inner product 
$\langle \cdot , \cdot \rangle: \mathcal{S}^N \times \mathcal{S}^N 
\rightarrow \R$, and $\mathcal{S}^N_+ \subset \mathcal{S}^N$ be the 
cone of positive semidefinite matrices. 
 
In several application domains, the SDPs to be solved involve 
\emph{discrete symmetries}.  In other words, there exists a finite 
group $G$, and an associated linear representation $\sigma:G 
\rightarrow \Aut(\mathcal{S}^N)$ for which certain relations hold.  The 
first requirement is that $\mathcal{S}^N_+$ maps onto itself under the 
group action, i.e., 
\begin{equation} 
\sigma(g)(\mathcal{S}^N_+)\subseteq \mathcal{S}^N_+, \quad\forall\,\, g\in G. 
\label{eq:SmapstoS} 
\end{equation} 
In subsequent sections we will require further specific properties for 
the representation $\sigma$, such as being induced by a representation 
$\rho:G \rightarrow GL(\R^N)$, via 
\begin{equation} 
\label{eq:induced} 
\sigma(g)(X):= \rho(g)^T X\rho(g),\quad \mbox{ for } X\in \mathcal{S}^N, \,\,g\in G, 
\end{equation} 
as this will be the case in most practical instances. Obviously, the 
required property (\ref{eq:SmapstoS}) is satisfied for such induced 
actions.  However, the results in the present section are valid for 
arbitrary linear representations $\sigma$, not necessarily of the form 
(\ref{eq:induced}) above. 
 
In order to describe and investigate our invariant semidefinite 
programs we need a few more definitions. A subset 
$\mathcal{L}\subseteq \mathcal{S}^N$ is called {\em invariant} with 
respect to $\sigma$ if $X\in \mathcal{L}$ implies $\sigma(g)(X)\in 
\mathcal{L},\,\forall\, g\in G$. We already required the property that 
the cone of positive definite matrices $\mathcal{S}^N_+$ is invariant. 
For a semidefinite program (\ref{eq:sdp1}) to be invariant we will 
demand that the set of feasible matrices be invariant. In other words, 
if a point belongs to the feasible set, so do its images under the 
group action: 
\[ 
X\in \mathcal{L}\cap \mathcal{S}^N_+ \quad \Rightarrow \quad 
\sigma(g)(X)\in \mathcal{L}\cap \mathcal{S}^N_+ \quad \forall\,g\in G. 
\] 
 
A linear mapping $F:\mathcal{S}^N\rightarrow \R$ is called {\em 
invariant} with respect to $\sigma$ if for all $X\in \mathcal{S}^N$ we 
have $F(X)=F\left(\sigma(g)(X)\right)$ for all $g\in G$. For the cost 
function we may restrict this property to the feasible set 
\[ 
\langle C , X \rangle   =   \bigl\langle C , \sigma(g)( X) \bigl\rangle,  
\qquad \forall\, g \in G, \,\,\forall\, X \in \mathcal{L}\cap \mathcal{S}^N_+. 
\] 
Equivalently, all points in a group orbit  
${\mathcal O}_X=\{ Y\in\mathcal{S}^N\,\,|\,\, \exists\, g \in G  \mbox{ with } \sigma(g)(X)=Y\}$  
have the same cost.

\begin{definition} \label{def:invSDP} 
  Given a finite group $G$, and associated linear representation  
$\sigma:G\rightarrow \Aut(\mathcal{S}^N)$, a semidefinite optimization problem 
of the form (\ref{eq:sdp1}) is called \emph{invariant} with respect to $\sigma$, if 
the following conditions are satisfied. 
\begin{enumerate} 
\item[(i)] The set of feasible matrices $\mathcal{L}\cap \mathcal{S}_+^N$ is   
invariant with respect to $\sigma$.  
\item[(ii)] The cost function $F(X)=\langle C , X \rangle$ is invariant with respect to $\sigma$.  
\end{enumerate} 
\end{definition} 
A related class of SDPs with symmetry has been defined independently
by Kanno \emph{et al.} in \cite{Murota}, where the authors study the
invariance properties of the search directions of primal-dual interior
point methods. Our results on the computational simplications that are
possible for these problems, however, are new as far as we can tell.

The group invariance of both the feasible set and the objective 
function suggests the possibility of considering a new optimization 
problem.  Concretely, for the invariant semidefinite program defined 
in Definition~\ref{def:invSDP}, an \emph{equivalent} optimization 
problem with a smaller number of variables and the same optimal 
solution can be obtained. The result is formally stated in 
Theorem~\ref{thm:restricted} below. To this end, in addition to 
invariant sets we also define {\em invariant} matrices as those that 
are fixed under the group action, i.e., they satisfy: 
\[ 
X=\sigma(g)(X), \qquad \forall\, g\in G. 
\] 
 
\begin{remark}  
  The linear representations $\sigma$ we are interested in are in 
  general orthogonal, i.e. the operators $\sigma(g) \in 
  \Aut(\mathcal{S}^N)$ satisfy that their inverse is equal to the 
  adjoint map, for all $g \in G$.  Then the matrix $C$ defining a 
  $\sigma$-invariant linear cost function $F(X)=\langle C,X\rangle$ is 
  itself an invariant matrix since 
\[  
\langle C, X\rangle= \bigl\langle C,\sigma(g) X\bigr\rangle =  
\bigl\langle \sigma(g)^T C, X\bigr\rangle= \bigl\langle \sigma(g^{-1}) C,X\bigr\rangle  
\quad \mbox{ for all }  X\in\mathcal{S}^n, g\in G. 
\] 
\label{rem:cmatrix} 
\end{remark}

We define the \emph{fixed-point subspace} of $\mathcal{S}^n$ (also 
called \emph{$G$-stable} or \emph{invariant}) as the subspace of all 
invariant matrices: 
\begin{equation} 
\mathcal{F}:= \{ X \in \mathcal{S}^n\,|\,\, X = \sigma(g) (X) \quad \forall\, g \in G\}, 
\label{eq:fixpointsub} 
\end{equation} 
and the associated semidefinite program: 
\begin{equation} 
F_{\sigma}:= \min_{X \in \mathcal{F}\cap \mathcal{L} \cap \mathcal{S}^n_+} \langle C , X \rangle, 
\label{eq:sdp2} 
\end{equation} 
obtained by restricting the feasible set of~(\ref{eq:sdp1}) to the 
fixed-point subspace~$\mathcal{F}$.

\begin{theorem} Given an orthogonal linear representation $\sigma:G\rightarrow \Aut(\mathcal{S}^n)$  
  of a finite group $G$ consider a semidefinite program which is 
  invariant with respect to $\sigma$. Then the optimal values $F$, 
  $F_\sigma$ of the original semidefinite program (\ref{eq:sdp1}) and 
  the fixed-point restricted semidefinite program (\ref{eq:sdp2}) are 
  equal. 
\label{thm:restricted} 
\end{theorem} 
 
\noindent 
{\bf Proof:} It is clear that $F \leq F_\sigma$, since every feasible 
matrix of (\ref{eq:sdp2}) is also a feasible matrix for~(\ref{eq:sdp1}). 
For the other direction, given any $\epsilon \geq 0$, consider a feasible 
matrix $X_1$ for (\ref{eq:sdp1}), such that $\langle C,X_1 \rangle = F 
+ \epsilon$. We will show that there exists a $X_\sigma\in\mathcal{F}\cap  
\mathcal{L}\cap\mathcal{S}_+^n$ with the same cost  
$\langle C,X_\sigma\rangle = F+ \epsilon$. This implies $F_\sigma\leq F+\varepsilon$. 
Define the ``group average'' (or Reynolds operator): 
\[ 
X_\sigma:= \frac{1}{|G|} \sum_{g\in G} \sigma(g) (X_1). 
\] 
We show first that $X_\sigma$ is a feasible point for (\ref{eq:sdp2}). 
Since $X_1\in \mathcal{L}\cap \mathcal{S}_+^n$ 
also $\sigma(g)(X_1)\in \mathcal{L}\cap \mathcal{S}_+^n$ by condition (i). 
First this means $\sigma(g)(X_1)\in \mathcal{L}$ and thus $X_\sigma\in \mathcal{L}$ 
because $\mathcal{L}$ is a linear space. 
Secondly, 
$\sigma(g)(X_1)$ is positive definite. Since the positive matrices form a convex cone 
each convex combination such as the group average  
defines another positive definite matrix.  
Thus $X_\sigma\in \mathcal{L}\cap \mathcal{S}_+^n$.  
Also, by 
construction, $X_\sigma \in \mathcal{F}$.  Finally, the cost function has the same value 
for $X_1$ and  $X_\sigma$ by condition~(ii):  
\[ 
\langle C, X_\sigma \rangle =  
\frac{1}{|G|} \sum_{g\in G} \bigl \langle C, \sigma(g)(X_1) \bigr \rangle =  
\frac{1}{|G|} \sum_{g\in G} \langle C, X_1 \rangle =  
\langle C, X_1 \rangle = F + \epsilon, 
\] 
so $F_\sigma \leq F + \epsilon$, from where the result follows. 
\hfill$\square$ 
 
The crucial property that will enable an efficient use of symmetry 
reduction in semidefinite programs is their \emph{convexity}. For 
general non-convex optimization problems, symmetries imply the 
existence of an orbit of minimizers (so there are still numerical 
advantages, see \cite{Wor94}), although the solutions themselves do 
not necessarily have any symmetric structure.  When convexity is 
present, as we have seen, the situation is much better, and the search 
for solutions can always be restricted \emph{a priori} to the 
fixed-point subspace.

\begin{example} 
The following very simple example illustrates the ideas. Consider 
a semidefinite program for $n=2$ with $\mathcal{L}=\mathcal{S}^2$  
and write $X=(X_{ij})$ for positive semidefinite matrices. Then 
\[ 
\min_{X\in \mathcal{S}_+^2} \langle C,X\rangle = 
\min_{X\in \mathcal{S}_+^2} (X_{11}+X_{22}) ,\quad \mbox{where} \quad 
C=  
\left[\begin{array}{cc} 
1 & 0 \\ 0 & 1  
\end{array}\right]. 
\] 
This semidefinite program is invariant under the group $C_2=\{id,s\}$ 
acting as $\sigma:C_2\rightarrow \Aut(\mathcal{S}^2)$ with 
\[ 
\sigma(s) 
\left[\begin{array}{rr} 
X_{11} & X_{12} \\ X_{12} & X_{22}  
\end{array}\right]= 
\left[\begin{array}{rr} 
X_{22} & - X_{12} \\ - X_{12} & X_{11}  
\end{array}\right].  
\] 
The fixed point subspace is  
$\mathcal{F}= \{X\in\mathcal{S}^2\,\,|\,\,X_{22}=X_{11},X_{12}=0\}$.   
By Theorem~\ref{thm:restricted} the original problem is equivalent to the 
trivial \emph{linear} program: 
\[ 
\min_{X\in \mathcal{F}\cap \mathcal{S}_+^2}\langle C,X\rangle= 
\min_{X_{11}\geq 0}(2 X_{11}). 
\]\hfill $\diamondsuit$ 
\end{example} 
For this particular example, the same conclusion could obviously be 
obtained by using elementary arguments. As we will see next, these 
basic techniques, in combination with the algebraic properties of the 
group representations that usually arise in practice, can be 
successfully applied to much more complicated examples.

\section{Induced actions and block diagonalization} 
\label{sec:blk} 

In most practical cases, the representation $\sigma$ of $G$ on 
$\mathcal{S}^N$ has a quite special form, since it arises as a 
consequence of a group action of $G$ on $\R^N$. This additional 
structure enables a computationally convenient characterization of the 
fixed point subspace. Essentially, under a convenient change of 
coordinates, every invariant matrix will be block diagonal. The 
details are provided below. 
 
The semidefinite programs we investigate in this section are invariant with respect to 
a linear representation $\sigma:G \rightarrow \Aut(\mathcal{S}^N)$ which 
is \emph{induced} by an action $\rho:G \rightarrow GL(\R^N)$, via: 
\begin{equation} 
\sigma(g) (X) := \rho(g)^T X \rho(g), \quad \forall\, g\in G. 
\label{eq:indrep} 
\end{equation} 
We assume that the representation $\rho$ is \emph{orthogonal}, that 
is, $\rho(g)^T \rho(g) = \rho(g) \rho(g)^T = I_N$, for all $g \in 
G$. Then, the restriction to the fixed point subspace 
(\ref{eq:fixpointsub}) takes the form 
\begin{equation} 
\rho(g) X  - X \rho(g) =0, \quad \forall\,\, g \in G. 
\label{eq:commute} 
\end{equation} 
Matrices commuting with the action of a group as in (\ref{eq:commute}) 
occur in many situations in mathematics and applications.  Using the 
Schur lemma from representation theory, it can be shown that all such 
matrices can be put into a \emph{block diagonal} form by a linear 
change of coordinates. In the following we briefly recall this 
well-known result and discuss its consequences for our invariant 
semidefinite programs.

In what follows, we consider first complex representations, as the theory is 
considerably simpler in this case.  For a finite group $G$ there are 
only finitely many inequivalent irreducible representations 
$\vartheta_1 , \ldots,\vartheta_h $ of dimensions (or degree) 
$n_1,\ldots,n_h$, respectively. The degrees $n_i$ divide the group 
order $|G|$, and satisfy the constraint $\sum_{i=1}^h n_i^2 = |G|$. 
Every linear representation of $G$ has a canonical decomposition 
\cite{Serre,Fassler} as a direct sum of irreducible representations, 
i.e.: 
\[ 
\rho = m_1 \vartheta_1 \oplus m_2 \vartheta_2 \oplus \cdots \oplus m_h \vartheta_h,  
\] 
where $m_1, \ldots, m_h$ are the \emph{multiplicities}. According to this there is 
an \emph{isotypic decomposition} 
\[ 
\C^n=V_1\oplus \cdots \oplus V_h, \quad V_i=V_{i1}\oplus \cdots \oplus V_{in_i}. 
\] 
A basis of this decomposition transforming with respect to the
matrices $\vartheta_i(g)$ is called \emph{symmetry adapted} and can be
computed using the algorithm presented in \cite{Serre,Fassler}. This
basis defines a change of coordinates by a matrix $T$ collecting the
basis as columns.  By Schur's lemma $T^{-1} X T$ has block diagonal
form with one block $X_i$ for each irreducible representation of
dimension $m_in_i$ which further decomposes into $n_i$ equal blocks
$B_i$ of dimension~$m_i$. that is:
\[ 
T^{-1}XT =  
\left[\begin{array}{ccc} 
X_1 & & {\bf 0}\\ 
& \ddots &\\ 
{\bf 0} && X_h 
\end{array}\right], \qquad X_i=\left[\begin{array}{ccc} 
B_i &&{\bf 0}\\ 
&\ddots\\ 
{\bf 0}&& B_i 
\end{array}\right]. 
\] 
 
For our applications of semidefinite programs, the problems are
usually presented in terms of real matrices, and therefore we would
like to use real coordinate transformations. In fact a generalization
of the classical theory to the real case is presented in
\cite[p. 106--109]{Serre}, and has been used many times in bifurcation
theory \cite{Golubitsky}).  If all $\vartheta_i(g)$ are real matrices
the irreducible representation is called {\em absolutely irreducible}.
Otherwise, for each $\vartheta_i$ with complex character its complex
conjugate will also appear in the canonical decomposition.  Since
$\rho$ is real both will have the same multiplicity and real bases of
$V_i+\bar V_i$ can be constructed. So two complex conjugate
irreducible representations form one real irreducible representation
of {\em complex type}.  There is a third case, real irreducible
representations of {\em quaternonian type}, seldom seen in practical
examples. A real irreducible representation of quaternonian type
decomposes over the complex numbers as twice a complex irreducible
representation with real character, but complex representation
matrices.
 
Since $\rho$ is assumed to be orthogonal the matrix $T$ can also be 
chosen to be orthogonal, too. Thus symmetric matrices $A$ are transformed 
to symmetric matrices~$T^TAT$. For symmetric matrices the block 
corresponding to a representation of complex type or quaternonian type
simplifies to a collection of equal subblocks. 
 
The change of coordinates affects both the explicit description of the 
affine linear space $\mathcal L$ and the cost function $F(X)=\langle 
C,X\rangle$. 
 
\begin{theorem} 
  Let $G$ be a finite group with $h$ real irreducible representations. 
  Let $\sigma$ 
  be a linear representation of $G$ on $\mathcal{S}^N$ induced as in 
  (\ref{eq:indrep}) by an orthogonal linear representation $\rho$ on 
  $\R^N$.  Given a semidefinite program which is invariant with 
  respect to $\sigma$ then there exists an explicitly computable 
  coordinate transformation under which all invariant matrices are 
  block diagonal. As a consequence, the semidefinite program can be 
  simplified to 
\begin{equation} 
\min_{X\in \mathcal{L}, X_i \in \mathcal{S}^{m_i}_+}  
\sum_{i=1}^h n_i\langle C_i , X_i \rangle, 
\label{eq:sdp3} 
\end{equation} 
where $m_i$  
are the multiplicities in $\rho$ 
and $n_i$ the dimensions of the real irreducible representations. 
$C_i$ denote the blocks of $C$ as well as $X_i$  
denote the blocks of $X$ in the new coordinates, respectively.   
\label{thm:newsdp} 
\end{theorem} 
The semidefinite program in Theorem~\ref{thm:newsdp} is simply the fixed-point 
restricted semidefinite program (\ref{eq:sdp2}) rewritten in a convenient coordinate 
frame in which the test of positive definiteness can be done much more 
efficiently.  
 
\begin{remark} 
In solving the semidefinite programs, we only need to consider 
\emph{one block} from \emph{each isotypic component}. In other words, 
instead of a large SDP of dimension $\sum_i n_i m_i$, we solve instead 
$h$ coupled semidefinite programs of dimension $m_i$.  Note also that 
some, or many, of the multiplicities $m_i$ can be zero. 
\end{remark} 

\begin{example} 
This example illustrates both procedures  
(Theorems \ref{thm:restricted} and \ref{thm:newsdp}).  
Consider the SDP with problem data given by: 
\[ 
C=\left[\begin{array}{ccc}0&0&0\\0&1&0\\0&0&1\end{array}\right] 
\quad  \mbox{and} \quad 
\mathcal{L}=\left\{X\in\mathcal{S}^3\,\,|\,\, \exists \, a,b,c_1,c_2,d\in\R  \mbox{ with } 
X=\left[\begin{array}{ccc} 
a & b & b \\ b & c_1 & d \\ b & d & c_2   
\end{array}\right] \right\}. 
\] 
Then the semidefinite program 
\begin{equation}\label{eq34} 
\min_{X\in \mathcal{L}\cap\mathcal{S}_+^3}\langle C,X\rangle= \min_{X\in \mathcal{L}\cap\mathcal{S}_+^3} (c_1 + c_2)  
\end{equation} 
is invariant under the action of the group $S_2$ given by simultaneous 
permutation of the last two rows and columns. To restrict the problem 
to the fixed point subspace, we can impose the constraint $c_1 = c_2 = c$, 
obtaining: 
\begin{equation}\label{eq35} 
\min_{X\in \mathcal{F} \cap \mathcal{L}\cap \mathcal{S}_+^3} (2 c)  
\quad \mbox{where}\quad 
\mathcal{F}=\left\{X\in \mathcal{S}^3\,\,|\, 
\exists \, a,b,c,d\in\R  \mbox{ with } 
\,X= 
\left[\begin{array}{ccc} 
a & b & b \\[-0.2cm] b & c & d \\[-0.2cm] b & d & c   
\end{array}\right]\right\}. 
\end{equation} 
$S_2$ has two irreducible representations of dimension $n_1=n_2=1$ 
which are also real irreducible. The multiplicities in the current representation are 
$m_1=2,\,m_2=1$. 
A symmetry adapted orthogonal coordinate transformation is given by 
the matrix 
\[ 
T=  
\left[\begin{array}{rrr} 
1 & 0 & 0 \\[-0.2cm] 0 & \alpha & \alpha \\[-0.2cm] 0 & \alpha & -\alpha   
\end{array}\right],  
\qquad  \alpha = \frac{1}{\sqrt{2}}. 
\] 
Now, the block diagonalization procedure can be applied, and the 
constraints are simplified to: 
\[ 
T^TCT=\left[\begin{array}{cc}C_1&0\\0&C_2\end{array}\right]\quad 
\mbox{with} \quad 
C_1=\left[\begin{array}{cc} 
0 & 0 \\[-0.2cm]  0  & 1 
\end{array}\right]  \quad  
\mbox{and} \quad  
C_2=1 \quad \mbox{and}  
\] 
\[ 
T^TXT=\left[\begin{array}{cc}X_1&0\\0&X_2\end{array}\right]\quad 
\mbox{with} \quad 
X_1=\left[\begin{array}{cc} 
a & \sqrt{2} b \\ \sqrt{2} b & c+d  
\end{array}\right]  \quad  
\mbox{and} \quad  
X_2=c-d.  
\] 
By Theorem~\ref{thm:newsdp} the optimization problems (\ref{eq34}) and (\ref{eq35}) are  
equivalent to  
\begin{eqnarray} 
\min_{X_1\in\mathcal{S}_+^2,X_2\geq 0} 
\bigl(\langle C_1,X_1\rangle+\langle C_2,X_2\rangle\bigr)& = & 
\min_{X_1\in\mathcal{S}_+^2,X_2\geq 0}\bigl((c+d)+(c-d)\bigr)\nonumber\\ 
&=& 
\min_{X_1\in\mathcal{S}_+^2,X_2\geq 0} (2 c). 
\end{eqnarray} 
\hfill$\diamondsuit$ 
\end{example}

\section{Sum of squares for invariant polynomials} 
\label{sec:polyinvar} 

The third level of specialization occurs when the semidefinite programs 
to be solved correspond to the computation of a sum of squares 
decomposition of multivariate polynomials \cite{pablo}. In this case, 
certain aspects of algebraic invariant theory 
\cite{Sta79,sturmfelsinvariant,DerKem02} appear very naturally in the picture. 
 
Assume we are interested in finding a sum of squares decomposition of
a polynomial $f(\mathbf{x})$ of degree $2d$ in $n$ variables which is
invariant with respect to a linear representation
$\vartheta:G\rightarrow GL(\R^n)$, i.e.,
$f(\mathbf{x})=f\left(\vartheta(g)\mathbf{x}\right),\,\, \forall\,
g\in G$. The set of all such invariant polynomials is the
\emph{invariant ring} $\R[\mathbf{x}]^G$. The representation
$\vartheta$ on $\R^n$ induces another linear representation $\tau$ of
$G$ on the polynomial ring $\R[\mathbf{x}]$ by
$\tau(g)\bigl(p(\mathbf{x})\bigr)=p\bigl(\vartheta(g^{-1})\mathbf{x}\bigr)$.
In this section, we study how this structure can be exploited
computationally.
 
\begin{example} 
On a first thought, one could jump to the conclusion that since the 
polynomial $f$ is invariant under the group action, we could further 
restrict the polynomials $p_i$ in the representation 
(\ref{eq:sosdecomp}) to being invariant, i.e., those in the trivial 
isotypic component. But this is not true, as the following simple 
example illustrates.  Consider the polynomial: 
\[  
f(x):= x^2 + (x-x^3)^2 = x^6 -2 x^4 + 2 x^2,  
\] 
invariant under the $C_2$ action generated by $x \rightarrow -x$.  The 
invariant ring is generated by $\theta(x)=x^2$, so the term $x^6$ 
cannot occur as the leading term of a sum of squares of invariant 
polynomials since $\theta^2=x^4$ and $(\theta^2)^2=x^8$. 
\hfill$\diamondsuit$ 
\label{rem:z2} 
\end{example} 
 
Instead of restricting to invariant polynomials we need to deal with 
the full polynomial ring.  Analog to the isotypic decomposition of 
$\R^n$ in Section~\ref{sec:blk} the polynomial ring also has an 
isotypic decomposition 
\begin{equation}\label{eq:ringiso} 
\R[\mathbf{x}]=V_1\oplus V_2\oplus\cdots \oplus V_h, 
\quad V_i=V_{i1}\oplus\cdots\oplus V_{in_i}. 
\end{equation} 
Here the isotypic component $V_1$ with respect to 
the trivial irreducible representation is the invariant ring 
$\R[\mathbf{x}]^G$.  The elements of the other isotypic components are called 
semi-invariants \cite{Sta79}. 
 
In order to find a symmetry adapted basis of an isotypic component of 
$\R[\mathbf{x}]$ corresponding to a real irreducible representation $\vartheta_i$ 
of dimension $n_i\geq 2$ one would like to compute equivariants. 
 
A polynomial mapping $f:\R^n\rightarrow \R^{n_i}$ with $f_i\in\R[\mathbf{x}]$ and 
\[ 
f\left(\vartheta(g) \mathbf{x}\right)=\vartheta_i(g)f(\mathbf{x}),\quad \forall\, g\in G, 
\] 
is called $\vartheta$-$\vartheta_i$-{\em equivariant}. 
 
Then the components $f_1,\ldots,f_{n_i}$ form a basis which transforms 
with respect to $\vartheta_i$, i.e. they form a symmetry adapted basis 
with respect to $\vartheta_i$. Thus $f_1,\ldots,f_{n_i}\in V_i$ and 
$f_1\in V_{i1},\ldots,f_{n_i}\in V_{in_i}$. 
There is a comprehensive body of theory that describes the invariants 
and equivariants \cite{Sta79,gat}. 
 
The following lemma is the key of our developments: 
\begin{lemma}[\cite{Sta79,gat}] 
  The set of $\vartheta$-$\vartheta_i$-equivariant forms is a finitely-generated 
  \emph{module} over the invariant ring $\R[\mathbf{x}]^G$. 
\end{lemma} 
By this lemma each $V_{ij}$ is a finitely generated module over 
$\R[\mathbf{x}]^G$. The algorithm for the computation of a symmetry adapted 
basis of $\R^n$ specializes to the computation of invariants and 
equivariants by analogue projections as outlined in 
\cite{Sta79,Ga96,Wor94}. 
 
There is much more known about the algebraic structures of the 
isotypic components $V_i$ of the polynomial ring which also helps to 
find generators.  We will deal with this in Section~\ref{sec:invar}. 
Within this section we just mention that the dimensions of the 
corresponding vector spaces (and therefore, the reduced SDPs) can be 
easily computed. 
 
The natural grading of $\R[\mathbf{x}]$ induces a grading on the $V_i$. 
Each $V_i$ decomposes into vector spaces of homogeneous polynomials. 
\[ 
V_i=\bigoplus_{k=0}^\infty (V_i)_k 
\] 
The formal sum  
\[ 
\psi_i(\xi) := \sum_{k=0}^\infty \dim\bigl((V_i)_k\bigr)\xi^k 
\] 
in the unknown $\xi$ is the well-known Hilbert-Poincar\'{e} 
series.  For the invariant and semi-invariants it is usually called 
the \emph{Molien} series. It can be computed as: 
\begin{equation} 
\psi_i(\xi) =  
\frac{1}{|G|}  
\sum_{g \in G} \frac{\trace{\vartheta_i(g)}^* }{\det(I - \xi \vartheta(g))} 
\end{equation} 
Software is available for the computation of these series and the
generators of the invariant ring and the module of equivariants.
While the Molien series in available in Singular, Magma and Macaulay 2
\cite{M2} the more specialized Hilbert-Poincar\'{e} series for
semi-invariants can be computed with special Maple packages
\cite{GaGu96symmetry,invar} and in Magma.  That means that for the
isotypic components $(V_i)_k$ the dimension and a vector space basis
can be easily computed.
 
Recall from Section~\ref{sec:sossdp} that we are interested in finding 
for a given $\vartheta$-invariant polynomial $f$ a positive definite 
matrix $A\in \mathcal{S}^N$ such that 
\[ 
f(\mathbf{x})=Y^T A Y, 
\] 
where we took $\lambda=0$ without loss of generality. The 
representation $\vartheta$ induces the orthogonal representation 
$\rho$ on $\R[\mathbf{x}]_{\leq d}$ and thus $\sigma$ on $\mathcal{S}^N$ by 
(\ref{eq:indrep}).  By Theorem~\ref{thm:newsdp} there exists an 
orthogonal change of coordinates $T$ such that 
\[ 
f(\mathbf{x})=Y^TTT^T A TT^TY=(T^TY)^T(T^TAT)(T^TY), 
\] 
where $T^TY=(v_1,\ldots,v_h)$ with $v_i=(v_{i1},\ldots,v_{in_i})$  
with $v_{ij}\in (\R[\mathbf{x}])^{m_i}$. Each component of $v_{ij}$ is an
element of $V_{ij}$ in~(\ref{eq:ringiso}). 
Again Schur's lemma gives
\[ 
T^TAT=\left[\begin{array}{ccc} A_1 && \mathbf{0}\\ 
& \ddots&\\ 
\mathbf{0}&&A_h\end{array}\right] \quad 
A_i=\left[\begin{array}{ccc}Q_i&& \mathbf{0}\\ 
&\ddots&\\ 
\mathbf{0}&&Q_i\end{array}\right]. 
\] 
The representation of the SDP cost function $F(X)=\langle E_{11},X\rangle=f(0)-\lambda$ 
is not affected by this change of coordinates.

If $f$ can be represented as a sum of squares, then there exists a 
representation of the form: 
\begin{equation} 
f(\mathbf{x})=\sum_{i=1}^h \sum_{j=1}^{n_i} 
\bigl( v_{ij}(\mathbf{x})\bigr)^T Q_i v_{ij}(\mathbf{x}). 
\label{eq:form} 
\end{equation} 
Notice that in the expression above, each of the $h$ terms correspond 
to a different irreducible representation of the group $G$. Of course, 
just like in the non-symmetric case, the representation is in general 
not unique. 
\begin{theorem} 
  Let $f$ be a $\vartheta$-invariant polynomial.  Then, $f(x)-f^{sos}$ 
  can be expressed as a sum of squares of invariant and semi-invariant 
  polynomials. 
\label{thm:sosinv} 
\end{theorem} 
\begin{example}  
\label{ex:dihedral} 
In this example we analyze a variation of the Robinson form, described 
in \cite{Reznick}. This instance has an interesting dihedral symmetry, 
and it is given by: 
\begin{equation} 
f(x,y) = x^6 + y^6 - x^4 y^2 - y^4 x^2 - x^4 - y^4 - x^2 - y^2 + 3 x^2 y^2 + 1.   
\end{equation} 
This polynomial is a dehomogenization of Robinson's form, i.e., 
$f(x,y) = R(x,y,1)$. It is invariant with respect to a representation 
of the dihedral group $D_4$ of order $8$, with the actions on $\R^2$ generated by: 
\[  
\vartheta(d)(x,y) = (-y,x), \qquad  \vartheta(s)(x,y) = (y,x). 
\] 
The generators $d, s$ have period $4$ and $2$, respectively. They 
satisfy the commutation relationship $s = dsd$. The group $D_4$ has 
five irreducible representations, shown in Table~\ref{tab:irrepd4}, of 
degrees $1,1,1,1,$ and $2$ \cite{Serre,Fassler}.  All of them are 
absolutely irreducible. 
 
A \emph{symmetry-adapted basis} for the five isotypic components can 
now be obtained, using the algorithm described in \cite{Fassler}, 
obtaining the corresponding basis vectors: 
\begin{equation} 
\begin{array}{rclrcl} 
v_1 &=& ( 1,x^2+y^2)  \hspace*{2cm} &v_{51} &=& ( x,x^3,x y^2) \\ 
v_2 &=& \emptyset  & v_{52} &=& ( y,y^3,y x^2). \\ 
v_3 &=& ( xy )  \\ 
v_4 &=& ( x^2 - y^2 ) 
\end{array} 
\label{eq:d4basis} 
\end{equation} 
\begin{table} 
\begin{center} 
\begin{tabular}{c|c|c} 
  & $d$ & $s$ \\ \hline 
$\vartheta_1$ &  1 & 1 \\ 
$\vartheta_2$ &  1 & $-1$  \\ 
$\vartheta_3$ &  $-1$ & $ 1 $ \\ 
$\vartheta_4$ &  $-1$ & $-1$  \\ 
$\vartheta_5$ &  
$\left[\begin{array}{cc} 
0 & -1 \\ 1 & 0
\end{array}\right] $ 
& 
$\left[\begin{array}{cc} 
0 & 1 \\ 1 & 0 
\end{array}\right]$  
\end{tabular} 
\end{center} 
\caption{Irreducible unitary representations of the dihedral group $D_4$.} 
\label{tab:irrepd4} 
\end{table} 
The corresponding Molien series for each irreducible component all 
have the common denominator $q(\xi) = (1-\xi^2) 
(1-\xi^4)$, and are given by: 
\begin{equation} 
\psi_1(\xi)  = \frac{1}{q(\xi)}, \quad 
\psi_2(\xi) = \frac{\xi^4}{q(\xi)}, \quad 
\psi_3(\xi) = \frac{\xi^2}{q(\xi)}, \quad 
\psi_4(\xi) = \frac{\xi^2}{q(\xi)}, \quad 
\psi_5(\xi) = \frac{\xi+\xi^3}{q(\xi)}. 
\end{equation} 
Notice that from these, by expanding the series, we can directly read 
the dimensions of the bases (\ref{eq:d4basis}) (and therefore the 
sizes of the corresponding SDPs), but much more important, we can do 
this for \emph{any} other polynomial with the same symmetry group. As 
expected, $\sum_i n_i \cdot \psi_i(\xi) = 
\frac{1}{(1-\xi)^2}$, the Hilbert series of the ring of 
polynomials in two variables. 
  
As suggested by Theorem~\ref{thm:newsdp}, after this symmetry 
reduction, the resulting semidefinite programs are much simpler. The 
affine subspace $\mathcal{L}$ is given by all matrices $Q\in 
\mathcal{S}^{10}$ which have block diagonal form with blocks 
$Q_1,Q_3,Q_4,Q_5,Q_5$ restricted to be 
\[ 
\begin{array}{rclrcl} 
Q_1& = & \left[ 
\begin{array}{cc} 
1-\lambda & c_1 \\ c_1 & c_2 
\end{array}\right], \hspace*{2.5cm} &  
Q_5& = & \left[\begin{array}{ccc} 
-1-2 c_1 & c_3 & c_4 \\ c_3 & 1 & c_5 \\ c_4 & c_5 & -1 - 2 c_5 
\end{array}\right],  
\\[0.3cm] 
Q_3 & = & 1 - 4 c_2 -4 c_3 -4 c_4,   & 
Q_4 & = & - 1 - c_2 -2 c_3,     
\end{array} 
\] 
where $c_1,c_2,c_3,c_4,c_5$ are real numbers. 
Then the cost function $1-\lambda$ is minimized over the set of 
feasible matrices $\mathcal{L}\cap\mathcal{S}_+^{10}$. 
That means each block $Q_i$ is restricted to be a positive semidefinite matrix. 
 
The optimal solution has $\lambda^* = -\frac{3825}{4096}$, with the 
corresponding matrices being
\begin{equation} 
Q_1=\left[ 
\begin{array}{rr} 
\frac{7921}{4096} & -\frac{89}{128} \\[0.2cm] -\frac{89}{128} & \frac{1}{4} 
\end{array}\right]=L_1L_1^T, \quad Q_3=Q_4=0, \quad  
Q_5=\left[\begin{array}{ccc} 
\frac{25}{64} & -\frac{5}{8} & \frac{5}{8} \\[0.1cm] 
-\frac{5}{8} & 1 & -1 \\[0.1cm] 
\frac{5}{8} & -1 & 1 \\  
\end{array}\right]=L_5L_5^T, 
\label{eq:d4matrices}  
\end{equation} 
with Cholesky factors 
\[ 
L_1=\left[\begin{array}{cc} 0 & -\frac{89}{64}\\[0.2cm] 0& \frac{1}{2} 
\end{array}\right], \quad  
L_5=\left[\begin{array}{ccc} 0 & 0& -\frac{5}{8} \\0 & 0&1\\ 
0 &0&-1\end{array}\right]. 
\] 
The previous theorem then yields that the sum of squares is 
\[ 
f(x,y) - \lambda^* =  
\left(-\frac{89}{64}+\frac{1}{2}(x^2+y^2)\right)^2 +  
\left(-\frac{5}{8} x+x^3-y^2 x\right)^2 +  
\left(-\frac{5}{8} y+y^3-x^2y\right)^2.  
\] 
Notice that the last two terms are \emph{not} invariant under the 
group action. However, their sum is, and this property can be fully 
exploited, as we shall see in Section~\ref{sec:invar}.\hfill$\diamondsuit$ 
\end{example}

\begin{example} 
  Consider the biquadratic form in six variables: 
\begin{eqnarray} 
B({\bf x};{\bf y}) &= & 
x_1^2 y_1^2 + x_2^2 y_2^2 + x_3^2 y_3^2 -  
2 \, (x_1 x_2 y_1 y_2+ x_2 x_3 y_2 y_3+x_3 x_1 y_3 y_1) +  \nonumber \\ 
&& + \, (x_1^2 y_2^2+x_2^2 y_3^2+x_3^2 y_1^2). 
\label{eq:choiform} 
\end{eqnarray} 
which is a minor modification of the celebrated example by Choi 
\cite{choiform}. Choi's original example was $B({\bf x},{\bf y}) + (x_1^2 
y_2^2+x_2^2 y_3^2+x_3^2 y_1^2)$, this modified version was studied by 
Choi and Lam \cite{ChoiLam}.  This form has several interesting 
properties, such as being nonnegative, even though it cannot be 
written as a sum of squares of polynomials. The associated matrix map 
$\mathcal{B}:\mathcal{S}^3 \rightarrow \mathcal{S}^3$ provided one of 
the first examples of a non-decomposable positive linear map. 
 
The form $B$ has numerous symmetries. The associated group has $96$ 
elements, with the corresponding generators given by: 
\begin{eqnarray*} 
(x_1,x_2,x_3,y_1,y_2,y_3) &\mapsto &(-x_1,x_2,x_3,-y_1,y_2,y_3)\\ 
\cdot\hspace*{2cm} &\mapsto&(x_1,-x_2,x_3,y_1,-y_2,y_3)\\ 
\cdot\hspace*{2cm}&\mapsto&(x_1,x_2,-x_3,y_1,y_2,-y_3)\\ 
\cdot\hspace*{2cm}&\mapsto&(x_1,x_2,x_3,-y_1,-y_2,-y_3)\\ 
\cdot\hspace*{2cm}&\mapsto&(x_3,x_1,x_2,y_3,y_1,y_2)\\ 
\cdot\hspace*{2cm}&\mapsto&(y_3,y_2,y_1,x_3,x_2,x_1). 
\end{eqnarray*} 
The group has $14$ irreducible representations: four one-dimensional, 
five of dimension two, four of dimension three, and one of dimension 
six ($4 \cdot 1^2 + 5 \cdot 2^2 + 4 \cdot 3^2 + 1\cdot 6^2 = 96$). The 
irreducible representations can be easily obtained using GAP4 
\cite{GAP4}. 
 
As mentioned, Choi showed that $B({\bf x};{\bf y})$ is nonnegative, but cannot be 
written as a sum of squares. Consider now the problem of checking 
whether $f({\bf x};{\bf y}) := (\sum_i x_i^2+y_i^2) B({\bf x};{\bf y})$ is a sum of squares. 
Clearly, all the symmetries of $B$ are inherited by $f$, since the 
additional factor is invariant under unitary transformations.  
 
We apply next the symmetry reduction procedure introduced in 
Section~\ref{sec:polyinvar}. After solving the SDPs we find a SOS 
decomposition where only two of the fourteen irreducible representations appear: a 
two-dimensional representation with multiplicity $m_1=3$ and a 
$6$-dimensional representation with multiplicity $m_2=4$. Indeed, $f$ 
can be written as: 
\begin{equation} 
f({\bf x};{\bf y}) =  \sum_{j=1}^2 v_{1j}^T Q_1 v_{1j} + \sum_{j=1}^6 v_{2j}^T Q_2 v_{2j}, 
\end{equation} 
where  
\[ 
\begin{array}{rclrcl} 
v_{11} &=& (x_3 y_1 y_2,x_2 y_1 y_3 ,x_1 y_2y_3 ),\hspace*{1.5cm} & 
v_{12} &=& (x_2 x_3 y_1 ,x_1 x_3 y_2,x_1 x_2 y_3 ),\\[0.2cm] 
v_{21} &=& (x_2 y_2y_3 ,x_3 {y_1 }^{2},x_1 y_1 y_3 ,x_3 {y_3 }^{2}),& 
v_{22} &=& (x_1 y_1 y_2,x_2 {y_3 }^{2},x_3 y_2y_3 ,x_2 {y_2}^{2}),\\[0.1cm] 
v_{23} &=& (x_1 x_2 y_2,{x_3 }^{2}y_1 ,x_1 x_3 y_3 ,{x_1 }^{2}y_1 ),& 
v_{24} &=& (x_1 x_3 y_1 ,{x_2 }^{2}y_3 ,x_2 x_3 y_2,{x_3 }^{2}y_3 ),\\[0.1cm] 
v_{25} &=& (x_2 x_3 y_3 ,{x_1 }^{2}y_2,x_1 x_2 y_1 ,{x_2 }^{2}y_2),& 
v_{26} &=& (x_3 y_1 y_3 ,x_1 {y_2}^{2},x_2 y_1 y_2,x_1 {y_1 }^{2}), 
\end{array} 
\] 
and $Q_1, Q_2$ are the positive semidefinite matrices 
\begin{equation} 
Q_1 = \left (\begin {array}{ccc} 1&-\half&-\half\\\noalign{\medskip}-\half&1&-\half 
\\\noalign{\medskip}-\half&-\half&1\end {array}\right ), \qquad 
Q_2 = \left (\begin {array}{cccc} 1&-1&\half&-\half\\\noalign{\medskip}-1&1&-\half 
&\half\\\noalign{\medskip}\half&-\half&1&-1\\\noalign{\medskip}-\half&\half&-1&1 
\end {array}\right ). 
\end{equation} 
\hfill$\diamondsuit$ 
\end{example}

\section{Computing with invariants}  
\label{sec:invar}  

The block-diagonalization techniques described in the previous section
directly enable a reduction in the computational burden of checking
SOS-based conditions. However, for efficiency reasons, in some cases
it is desirable to employ an alternative procedure where all
computations are done using only invariant elements, and if possible
never deal at all with the original variables.
  
The main benefit of computing only in the invariant ring is a much  
smaller number of algebraic operations, besides the fact that in many  
symmetric problems the original data is already expressed in a compact  
form, in terms of invariants (see Example \ref{ex:sottile}).  
Additionally, the resulting formulation has interesting theoretical  
properties, that shed more light on the nature of the sum of squares  
techniques. 
  
It is well-known that the invariant ring $\R[{\bf x}]^G$ is finitely 
generated as a $\R$-algebra.  Denoting a set of homogeneous generators 
by $g_1({\bf x}),\ldots,g_s({\bf x})$ this fact can be written as 
\[  
\R[{\bf x}]^G=\R[g_1,\ldots,g_s].  
\]  
Consequently, for each invariant polynomial $r({\bf x})$ there exists  
a polynomial $\tilde r(y_1,\ldots,y_s)$ such that  
\[  
r({\bf x})=\tilde r\bigl(g_1({\bf x}),\,\ldots\,,g_s({\bf x})\bigr).  
\]  
Given the polynomial generators $g_i$ the algorithmic computation of 
an appropriate $\tilde r$ is a standard task in computational algebra, 
for which several techniques are available. 
 
{\em Gr\"obner bases}: Besides special techniques for distinguished
algebra bases, the general {\em Gr\"obner bases} approach may be used
in this context.  In \cite{sturmfelsinvariant} a method is explained
based on slack variables $y_1,\ldots,y_s$ and a Gr\"obner basis of
$\langle y_1-g_1({\bf x}),\ldots, y_s-g_s({\bf x})\rangle$ with
respect to an elimination order. The division algorithm applied to $r$
gives a representation $\tilde r$.
  
There are two special sets of generators having advantageous properties.  
  
{\em SAGBI bases:} If the invariant ring admits a finite SAGBI basis 
\cite[Chapter 11]{Stbook}, $\tilde r$ is determined by the {\em 
  subduction algorithm}.  Since this technique does not use slack 
variables it is more efficient than the Gr\"obner basis approach.  A 
nice example is the case of the symmetric group, where the method is 
equivalent to the well-known algorithm for expressing an invariant 
polynomial in terms of the elementary symmetric functions.  However, 
this variant is not always possible, as there exist groups whose 
invariant rings have no finite SAGBI basis, for any term ordering.  We 
will not follow this approach here. 
  
{\em Hironaka decomposition:} It can be shown (see 
\cite{sturmfelsinvariant} and the references therein), that the 
invariant ring is Cohen-Macaulay, and has a Hironaka decomposition: 
\[  
\R[{\bf x}]^G = \bigoplus_{j=1}^t \eta_j({\bf x}) \, \R[\theta_1({\bf x}), \ldots, \theta_n({\bf x})],  
\]  
with $\theta_i({\bf x}), \eta_i({\bf x})$ being \emph{primary} and 
\emph{secondary invariants}, respectively. In particular the primary 
invariants are algebraically independent.  Every invariant polynomial 
can be written \emph{uniquely} in terms of $\theta_i,\eta_j$ with 
$\eta_j$ appearing linearly only. Algorithms for the computation of 
primary and secondary invariants are presented in 
\cite{sturmfelsinvariant,gat,DerKem02}. This unique representation 
may be found by using Gr\"obner bases with respect to special term 
orders or by {\em involutive bases} \cite{Plesken01}. 
 
In this paper we use the Hironaka decomposition (i.e., primary and 
secondary invariants) as the computational representation for the 
invariant ring. The main reasons are the advantages of having a unique 
representation, the available results on the structure on the isotypic 
components (e.g., Theorem~\ref{thm:modules} below), and the 
simplifications that occur in important special cases. 
 
For instance, in practice very often the group representations are
so-called \emph{reflection} groups.  This means that $\vartheta(G)$ is
generated by reflections with respect to hyperplanes. In this
situation the invariant ring is isomorphic to a polynomial ring (no
secondary invariants are needed) and the other isotypic components are
free modules over the invariant ring \cite{sturmfelsinvariant,Sta79}.
Our results are particularly appealing in this situation.

{\bf Example~\ref{ex:dihedral} continued:} {\em 
We revisit Example~\ref{ex:dihedral}, from this algebraic  
perspective. The action of $D_4$ is a reflection group and thus 
the invariant ring $S=\R[x,y]^{D_4}$ is generated by two primary  
invariants:  
\[  
\theta_1(x,y) = x^2+y^2 , \qquad \theta_2(x,y) = x^2 y^2. 
\]  
No secondary invariants are needed, therefore the invariant ring is 
isomorphic to a polynomial ring. Rewriting $f(x,y)$ in terms of the 
invariants using any of the methods described, we obtain $\tilde f\in 
\tilde T=\R[\theta_1,\theta_2]$ with: 
\begin{equation}  
f(x,y) = \tilde f(\theta_1(x,y),\theta_2(x,y)), \quad  \tilde f(\theta_1,\theta_2)=   
\theta_1^3 - \theta_1^2 - 4 \theta_1 \theta_2 - \theta_1 + 5 \theta_2 + 1.   
\end{equation}  
\hfill $\diamondsuit$ 
} 
 
A difficulty in trying to directly apply invariant-based techniques to 
sum of squares problems is that in a SOS decomposition the individual 
components (the squares) are \emph{not} necessarily invariant, e.g. 
Example~\ref{rem:z2}.  Nevertheless, the different terms corresponding 
to the same isotypic component can be conveniently ``bundled'' 
together, as we show next. 
 
The representation (\ref{eq:form}) may be written in a different way.  
Since $z^T Q z = \langle Q \,, z z^T\rangle$, then for each $i$:  
\[  
\sum_{j=1}^{n_i} \bigl(v_{ij}({\bf x})\bigr)^T Q_i \,v_{ij}({\bf x}) =  
\biggl\langle Q_i\,,\sum_{j=1}^{n_i} v_{ij}({\bf x})\bigl(v_{ij}({\bf 
  x})\bigr)^T\biggr\rangle  
=\bigl\langle Q_i, P_i({\bf x})\bigr \rangle   
\]  
using $P_i({\bf x}) = \sum_{j=1}^{n_i} v_{ij}({\bf x})(v_{ij}({\bf 
  x}))^T\in(\R[{\bf x}])^{m_i\times m_i}$.  Since the decomposition is 
symmetry-adapted, the entries 
$p^i_{kl}=\sum_{j=1}^{n_i}(v_{ij})_k(v_{ij})_l$ of the matrices 
$P_i=(p^i_{kl})_{kl}$ are scalar products of \emph{equivariant} 
polynomial vectors.  Therefore, because the representations are 
assumed to be orthogonal, the $p^i_{kl}({\bf x})$ are invariant 
polynomials, obtaining the representation: 
\[  
f({\bf x})=\sum_{i=1}^h \bigl\langle Q_i , P_i({\bf x}) \bigr\rangle, \qquad  
   P_i \in (\R[{\bf x}]^G)^{m_i \times m_i}.  
\]  
The matrices $P_i$ have components which are invariant polynomials, 
and therefore can be expressed in terms of primary and secondary 
invariants $\theta_i, \eta_i$. 
 
Since both $f$ and the entries of $P_i$ are invariant, there are
unique representations $\tilde f({\bf \theta},{\bf \eta})=\sum_{j=1}^t
\eta_j\tilde f_j({\bf \theta})$ and $\tilde p^i_{kl}({\bf \theta},{\bf
  \eta})$ being linear in $\eta$ such that
\[  
f({\bf x})=\tilde f\bigl(\theta({\bf x}),\eta({\bf x})\bigr)= 
  \sum_{j=1}^t \eta_j({\bf x})\,\tilde f_j\bigl(\theta_i({\bf x}),\ldots,\theta_n({\bf x})\bigr) 
, \quad  p^i_{kl}({\bf x})=\tilde p^i_{kl}\bigl(\theta({\bf x}),\eta({\bf x})\bigr). 
\] 
So the sum of squares of an invariant polynomial corresponds to a representation  
\[  
\tilde f({\bf \theta},{\bf \eta})=\sum_{i=1}^h \bigl\langle Q_i , \tilde P_i({\bf \theta},{\bf \eta}) \bigr\rangle,  
\qquad \tilde P_i \in (\tilde T)^{m_i \times m_i}, \quad  
\tilde T=\sum_{j=1}^t\eta_j\R[{\bf \theta}], 
\]  
where the matrices $Q_i$ are positive semidefinite. 
 
In order to compute in the coordinates $\theta_i,\eta_j$ we are left 
with the task to determine the matrices $\tilde P_i$.  To provide an 
efficient representation of these matrices, we will exploit the 
algebraic structure of the module of equivariants. 
\begin{theorem}(\cite{Wor94}) 
  For each irreducible representation $\vartheta_i$ the module $M_i$
  of $\vartheta$-$\vartheta_i$-equivariants is a free module over the
  ring in the primary invariants:
\[  
M_i =  \R[\theta]\cdot\{ b^i_1,\ldots,b^i_{r_i}\},  
\]  
where $b^i_\nu\in\R[{\bf x}]^{n_i}$ denote the elements in the module 
basis and $r_i$ is the rank of the module. 
\label{thm:modules} 
\end{theorem}  
 
For some special cases the ranks $r_i$ are known \emph{a priori}.  For 
instance, in \cite[Proposition 4.7]{Sta79} it is stated that for 
groups generated by pseudo-reflections, the module of equivariants 
associated with one-dimensional irreducible representations have 
$r_i=1$. 

In the case of a one-dimensional irreducible representation ($n_i=1$)
we have $M_i=V_i$. In general the spaces $V_{ij}$ are generated by the
components of the equivariants:
\[
V_{ij}=\R[\theta]\cdot\{(b_1^i)_j,\ldots,(b_{r_i}^i)_j\}.
\]
We choose the special vector space basis $v_i({\bf x})$ determined by the module basis, 
that is 
\begin{equation}\label{clever} 
\left[\begin{array}{c} 
v_{i1} \\  \vdots \\  v_{i n_{i}} 
\end{array}\right] =  
\left\{ 1, \theta_1, \ldots, \theta_n, \theta_1^2, \ldots \right\} 
\cdot 
\left\{ b_1^i , \ldots ,  b^i_{r_i} \right\}. 
\end{equation}

{\bf Example~\ref{ex:dihedral} continued:} {\em As stated before $D_4$ 
  has $5$ irreducible representations of dimensions $1$ and $2$. For 
  the one-dimensional representations the modules are generated by one 
  element: $r_2=r_3=r_4=1$.  The generators are $xy(x^2-y^2)$ and $xy$ 
  and $x^2-y^2$, respectively.  The fifth irreducible representation 
  has dimension $2$; the free module has two generators, namely 
\[ 
b^5_1(x,y)=\left(\begin{array}{c}x\\y\end{array}\right),\quad 
b^5_2(x,y)=\left(\begin{array}{c}x^3\\y^3\end{array}\right). 
\] 
To make explicit this module structure, we  
replace the basis elements $v_{51}, v_{52}$ in (\ref{eq:d4basis}) by the  
equivalent ones:  
\[  
v_{51} = \bigl(x,x^3,(x^2+y^2)x\bigr)= \bigl(x,x^3,\theta_1 x\bigr), \quad  
v_{52} = \bigl(y,y^3,(x^2+y^2)y\bigr )=\bigl(y,y^3,\theta_1 y\bigr ).   
\]  
\hfill $\diamondsuit$ 
}

Now we show how the information in the matrix $\tilde P_i$ can be 
``compressed'' to a much smaller one, of size equal to the module rank 
$r_i$. 
 
Let $\Pi_i$ be the $r_i \times r_i$ principal subblock of $\tilde P_i$
defined by the first $r_i$ rows and columns. Without loss of
generality we can assume that the size of the matrix $\tilde P_i$ is a
multiple of $r_i$ (otherwise we add rows and columns of zeros to
$Q_i$).  By the choice of vector basis in (\ref{clever}), the other
$r_i\times r_i$-subblocks of $\tilde P_i$ are just monomial multiples
of $\Pi_i$.  The monomials are just monomials in ${\bf \theta}$ only.
Similarly, we split $Q_i$ into $r_i\times r_i$-subblocks. We can then
rewrite the inner products as operations between the polynomial
matrices $\Pi_i$ and SOS matrices $S_i$, defined in terms of $Q_i$. We
formally state our result as follows:
 
\begin{theorem} 
Let $f(\mathbf{x})$ be a SOS polynomial, that is invariant under the 
action of a finite group $G$, and $\theta_i, \eta_i$ be primary and 
secondary invariants of the corresponding invariant ring. Then, 
$\tilde f(\theta,\eta)$ has a representation of the form: 
\begin{equation} 
\tilde f({\bf \theta},{\bf \eta}) =\sum_{i=1}^h\bigl\langle  S_i({\bf \theta} ),  
\Pi_i({\bf \theta},{\bf \eta})\bigr\rangle, 
\label{eq:matrixsosrep}
\end{equation}
where $\Pi_i \in \R[\theta,\eta]^{r_i \times r_i}$ depend only on the
 group action, and $S_i \in \R[\theta]^{r_i \times r_i}$ are SOS matrices. 
\label{thm:matrixsosrep} 
\end{theorem} 
{\bf Proof:} For notational simplicity, in this proof we drop the
index $i$ denoting the $i$-th irreducible representation.  For each
real irreducible representation the matrix $\Pi\in (\tilde T)^{r\times
  r}$ is defined in the following way. The entry
$\pi_{kl}(\theta,\eta)$ is the invariant polynomial $(b_k({\bf
  x}))^Tb_l({\bf x})$ represented in terms of fundamental invariants.
We represent the basis $v$ as in (\ref{clever}), where
$w=(w_1,\ldots,w_\nu)^T$ is the vector of monomials in $\theta$.
Using a weighted degree on $\R[\theta]$ and $\tilde T$ induced by the
degrees of $\theta_i({\bf x}), \eta_j({\bf x})$ it is clear that we
can take $w$ to contain all monomials in $\theta$ of degree less than
or equal to $\half(\deg(\tilde{f}) - \min_k \deg \pi_{kk} )$. In
particular, we take $w_1=1$.

As mentioned, in the chosen basis the matrix $\tilde P$ is formed by
blocks which are scalar copies of $\Pi$, having the structure $\tilde
P = w w^T \otimes \Pi$. As a consequence, we have:
\begin{eqnarray*}
\langle Q, \tilde P \rangle &=& 
\langle Q, w w^T \otimes \Pi \rangle = 
\langle \hat Q, \Pi \otimes w w^T  \rangle = 
\langle \hat Q, (I_r \otimes w) ( \Pi \otimes 1 ) ( I_r \otimes w)^T
\rangle \\
&=& 
\langle ( I_r \otimes w)^T \hat Q (I_r \otimes w),  \Pi \rangle = 
\langle S ,  \Pi \rangle.
\end{eqnarray*}
where $\hat Q$ is a permutation of the rows and columns of $Q$, and
$S$ is being defined by the last equality. It is clear that $S \in
\R[\theta]^{r \times r}$ is a SOS matrix as in
Definition~\ref{def:sosmatrix}, since it satisfies $\mathbf{y}^T
S\mathbf{y}= (\mathbf{y}\otimes w)^T \hat Q (\mathbf{y}\otimes w)$,
and $\hat Q$ is positive semidefinite since $Q$ is.
An alternative formula for $S$ can easily be obtained. Partition
$\hat Q=(\hat Q^{jl})_{j,l=1,\ldots,r}$ into $\nu\times \nu$-subblocks. 
This yields $(S)_{jl}=w^T \hat Q^{jl}w$ for $j=1,\ldots,r,l=1,\ldots,r$. 
\hfill\qed

When the numerical value of the invariants $\theta, \eta$ correspond
to \emph{real} values of the original variables, then the matrices
$\Pi_i(\theta,\eta)$ will always be positive semidefinite, since they
are constructed as sums of outer products (notice that they are
\emph{not} PSD for arbitrary values of $\theta, \eta$).  The $\Pi_i$
so generated therefore provide an ``easily parameterizable'' family of
nonnegative polynomials, constructed by contracting the matrices with
arbitrary SOS matrices $S_i$. Since inner products of positive
semidefinite matrices are always nonnegative, the representation
(\ref{eq:matrixsosrep}) provides a ``self-evident'' certificate of the
nonnegativity of $f$.

 \begin{remark}  
  The image of $\R^n$ under the invariants (the \emph{orbit space}) is  
  a semialgebraic set. Nonnegative symmetric polynomials naturally  
  correspond to polynomials that are nonnegative over this set; see  
  \cite{Procesi}.  The $\Pi_i$ are related to the defining  
  inequalities for this set, as well as its matching stratifications  
  (semialgebraic and by orbit-type). The representation obtained  
  closely matches what would be obtained by using the  
  Positivstellensatz, an issue that will be pursued elsewhere.  
\end{remark}  

There are a few practical details that help in the computation of
symmetric SOS representations as presented above. For instance, since
we know that the degree of the $(kk)$-entry of $S_i$ is bounded in
weighted degree by $\deg(f)-\deg(\pi^i_{kk})$, it is convenient to use
different monomial sets for the entries of $S$, rather than a common
set $w$ as described in the preceding proof. This is equivalent to
eliminating \emph{a priori} the zero rows and columns in $\hat Q$. Also,
using the unique representation property of the Hironaka
decomposition, the required identity (\ref{eq:matrixsosrep}) can be
expressed concisely as:
\[ 
\tilde f_j({\bf \theta}) =
\sum_{i=1}^h\bigl\langle  S_i({\bf \theta} ),  
\Pi_i^j({\bf \theta})\bigr\rangle, \qquad j = 1, \ldots, t,
\]
where $\Pi_i(\theta,\eta)=\sum_j\Pi_i^j(\theta)\,\eta_j$.

{\bf Example~\ref{ex:dihedral} continued:} {\em From the isotypic
  components, and using the new basis, we directly compute the
  matrices $\Pi_1,\Pi_2,\Pi_3,\Pi_4\in \R[\theta]^{1\times 1}$
  corresponding to the one-dimensional representations:
\[  
\Pi_1=  
1  
, \quad   
\Pi_2 = \theta_2 (\theta_1^2 - 4 \theta_2), \quad  
\Pi_3=\theta_2, \quad   
\Pi_4=\theta_1^2 - 4 \theta_2.   
\]  
For the fifth component we have with $w_5=(1,\theta_1)$ the extended basis: 
\begin{eqnarray*} 
v_{51} &=& \bigl(x,x^3,\theta_1 x,\theta_1 x^3\bigr) = (1,\theta_1)\otimes (x,x^3)  \\ 
v_{52} &=& \bigl(y,y^3,\theta_1 y, \theta_1 y^3\bigr) =  (1,\theta_1) \otimes (y,y^3),  
\end{eqnarray*} 
yielding 
\[  
\tilde P_5=\left[\begin{array}{cc|cc}  
\theta_1 & \theta_1^2 - 2 \theta_2 & \theta_1^2 & \theta_1(\theta_1^2 - 2 \theta_2)\\  
\theta_1^2-2 \theta_2 & \theta_1 (\theta_1^2 - 3 \theta_2) &   
\theta_1 (\theta_1^2-2 \theta_2) & \theta_1^2 (\theta_1^2 - 3 
\theta_2)\\  
\hline 
\theta_1^2 &   
\theta_1 (\theta_1^2-2 \theta_2) &  
\theta_1^3 & \theta_1^2(\theta_1^2 - 2 \theta_2)\\ 
\theta_1(\theta_1^2-2 \theta_2) & \theta_1^2 (\theta_1^2 - 3 \theta_2) &   
\theta_1^2 (\theta_1^2-2 \theta_2) & \theta_1^3 (\theta_1^2 - 3 
\theta_2)\\  
\end{array}\right]  =  
\left[\begin{array}{cc} 
\Pi_5 & \theta_1 \Pi_5 \\ 
\theta_1 \Pi_5 & \theta_1^2 \Pi_5  
\end{array}\right] 
.   
\]  
The subblock which appears multiplied repeatedly with $\theta_1$ and 
$\theta_1^2$ is 
\[  
\Pi_5=\left[\begin{array}{cc}  
\theta_1 & \theta_1^2 - 2 \theta_2 \\  
\theta_1^2-2 \theta_2 & \theta_1 (\theta_1^2 - 3 \theta_2) 
\end{array}\right].\]
Analogously, $Q_5, \hat Q_5$ have the subblock structure 
\[ 
Q_5=\left[\begin{array}{cc|cc} 
\frac{25}{64} & -\frac{5}{4} & \frac{5}{8}&0 \\[0.1cm] 
-\frac{5}{4} & 4 & -2 &0\\[0.1cm] 
\hline 
\frac{5}{8} & -2 & 1 &0\\  
0&0&0&0 
\end{array}\right], \qquad  
\hat Q_5=\left[\begin{array}{cc|cc} 
\frac{25}{64} & \frac{5}{8}& -\frac{5}{4} &0 \\[0.1cm] 
\frac{5}{8} &  1 &-2  &0\\  \hline
-\frac{5}{4} & -2 & 4  &0\\[0.1cm] 
0&0&0&0 
\end{array}\right]=
\left[\begin{array}{cc} 
\hat Q_{5}^{11} & \hat Q_5^{12} \\ 
(\hat Q_5^{12})^T & \hat Q_5^{22}  
\end{array}\right], 
\]
where the ordering in the latter corresponds to the basis
\[
\begin{array}{rcl} 
\tilde v_{51} &=& \bigl(x,\theta_1 x,x^3,\theta_1 x^3\bigr)=  (x,x^3)
\otimes (1,\theta_1),  \\ 
\tilde v_{52} &=& \bigl(y,\theta_1 y, y^3,\theta_1 y^3\bigr)=  (y,y^3)
\otimes (1,\theta_1). 
\end{array}
\]
We have the identity: 
\begin{eqnarray*} 
\langle Q_5, \tilde P_5 \rangle &=&  
\left\langle Q_5, 
\left[\begin{array}{cc}
1 & \theta_1 \\ \theta_1 & \theta_1^2 
\end{array}\right]
\otimes \Pi_5 \right \rangle = 
\left\langle \hat Q_5, 
\Pi_5 \otimes
\left[\begin{array}{cc}
1 & \theta_1 \\ \theta_1 & \theta_1^2 
\end{array}\right] 
\right \rangle \\
&=&
\langle
(I_2 \otimes w_5)^T \hat Q_5
(I_2 \otimes w_5), \Pi_5 \rangle 
= \langle S_5, \Pi_5 \rangle.
\end{eqnarray*} 
where the explicit expression for $S_5$ is given by:
\[
S_5=\left[ \begin{array}{cc}
w_5^T \hat Q^{5}_{11}w_5 & w_5^T \hat Q^{5}_{12}w_5\\
 w_5^T (\hat Q^{5}_{12})^Tw_5 & w_5^T \hat Q^{5}_{22}w_5
\end{array}\right].
\]
The construction guarantees $S_5$ is a SOS polynomial matrix, 
yielding the final decomposition:  
\[  
\tilde f(\theta_1,\theta_2) - \lambda^* =   
\left(-\frac{89}{64}+\frac{\theta_1}{2}\right)^2  \Pi_1 +   
\left \langle  
\left[\begin{array}{ccc}  
(\theta_1+\frac{5}{8})^2 & -\frac{5}{4}-2 \theta_1 \\  
-\frac{5}{4} - 2 \theta_1 & 4  
\end{array}\right] 
,  
\Pi_5  
\right \rangle.  
\]  
\hfill $\diamondsuit$ 
}

The next example illustrates the results in a case where secondary
invariants are involved.

\begin{example} 
Consider the cyclic group $C_4$ in four elements acting on $\R^2$ as: 
\[ 
(x,y) \mapsto (-y,x). 
\] 
The corresponding invariant ring $\R[x,y]^{C_4}$ is \emph{not} isomorphic 
to a polynomial ring in two variables. Instead, it is generated by  
primary and secondary invariants given by: 
\[ 
\theta_1 = x^2+y^2, \qquad 
\theta_2 = x^2 y^2, \qquad 
\eta_1 = 1, \qquad 
\eta_2 = x y (x^2 - y^2), 
\] 
that satisfy the relation (or \emph{syzygy}) $\eta_2^2 + 4 \theta_2^2-\theta_1^2 \theta_2 =0$. 
 
The group $C_4$ has four one-dimensional complex irreducible 
representations ($4 \cdot 1^2 = |C_4|$). Specializing to the field of real numbers, 
two of these ones merge, so the dimensions over $\R$ are 1,1, and 2. 
 
The isotypic components of $\R[x,y]$ are: 
\begin{eqnarray*} 
V_1 &=& \R[\theta_1,\theta_2] \cdot \{1,xy (x^2-y^2) \} \\ 
V_2 &=& \R[\theta_1,\theta_2] \cdot \{x y, x^2-y^2 \} \\ 
V_{3}&=& \R[\theta_1,\theta_2] \cdot  
\left\{  
\left(\begin{array}{c} x \\ y \end{array}\right) , 
\left(\begin{array}{c} x^3 \\ y^3 \end{array}\right)  
\right\}. 
\end{eqnarray*} 
The corresponding matrices $\Pi_i$ take the form: 
\[ 
\Pi_1 =  
\left[\begin{array}{cc}  
1 & \eta_2 \\  
\eta_2 & \theta_1^2 \theta_2 - 4 \theta_2^2 \end{array}\right], \qquad  
\Pi_2 =  
\left[\begin{array}{cc}  
\theta_2 & \eta_2 \\  
\eta_2 & \theta_1^2 - 4 \theta_2 \end{array}\right], \qquad  
\Pi_3 =  
\left[\begin{array}{cc}  
\theta_1 & \theta_1^2-2 \theta_2 \\  
\theta_1^2-2 \theta_2 & \theta_1^3 - 3 \theta_1 \theta_2 \end{array}\right]. 
\] 
Therefore, by Theorem~\ref{thm:matrixsosrep}, every $C_4$-invariant  SOS 
polynomial has a representation of the form 
\[ 
\tilde f(\theta_1,\theta_2,\eta_2) =  
\sum_{i=1}^3 \bigl\langle S_i(\theta_1,\theta_2),  \Pi_i(\theta_1,\theta_2,\eta_2) \bigr \rangle,  
\] 
where the $S_i$ are sum of squares matrices. 
\hfill$\diamondsuit$
\end{example} 
 
\begin{remark}  
  It is interesting to notice that some of the matrices $\Pi_i$ have
  simple interpretations. Consider the case of the symmetric group
  $S_n$. It is well-known that the lowest dimension where the
  alternating representation appears non-trivially is equal to
  $\binom{n}{2}$. In that case, the corresponding $\Pi$ is a scalar
  polynomial of degree $n(n-1)$, and equal to the discriminant
  $\prod_{i<j} (x_i-x_j)^2$ of the polynomial $\Pi_{i=1}^n
  (z-x_i)$. In fact, explicit formulas for the $\Pi$ corresponding to
  one-dimensional representations of reflective groups can be obtained
  explicitly; see \cite[Proposition 4.7]{Sta79}.
\end{remark}

\section{The algorithm}  

In the following we summarize the steps which are needed in order to 
efficiently compute a sum of squares representation of an invariant polynomial. 
Since several steps of the computation are useful for \emph{all} 
polynomials with a given invariance group, we distinguish between the 
operations that depend only on the specific group (which can be 
thought of as a preprocessing step), and those that operate on the 
concrete instance. 
 
{\bf Algorithm I} (Compute SOS generators) 
 
{\bf Input:} linear representation $\vartheta$ of a finite 
group $G$ on $\R^n$. 
 
\hspace*{0.7cm} {\bf 1.} determine all real irreducible representations of $G$.
 
\hspace*{0.7cm} {\bf 2.} compute primary and secondary invariants, $\theta_i,\eta_j$. 
 
\hspace*{0.7cm} {\bf 3.} for each non-trivial irreducible representation
compute the basis $b_1^i,\ldots,b_{r_i}^i$ \\
\hspace*{1.1cm} of the module of equivariants.  
 
\hspace*{0.7cm} {\bf 4.} for each irreducible representation $i$
compute the corresponding matrix $\Pi_i(\theta,\eta)$.

{\bf Output:} Primary and secondary invariants $\theta, \eta$. Matrices $\Pi_i$.

{\bf Algorithm II } (Sum of squares for invariant polynomials) 
 
{\bf Input:} Primary and secondary invariants $\theta, \eta$. Matrices $\Pi_i$ and $ f\in \R[{\bf x}]^G$. 
 
\hspace*{0.7cm} {\bf 1.} rewrite $f$ (if needed) in fundamental invariants 
giving $\tilde f(\theta,\eta)$. 
 
\hspace*{0.7cm} {\bf 2.} for each irr. repr. determine $w_i(\theta)$ and thus the structure of the matrices $S_i\in\R[\theta]$.  
 
\hspace*{0.7cm} {\bf 3.} find a feasible solution of the semidefinite
program corresponding to \\[0.2cm]
\hspace*{1.1cm} the constraints $\{S_i \mbox{
  are SOS matrices}, \,\,
\tilde f=\sum_{j=1}^h\langle S_i, \Pi_i\rangle\}$.
 
{\bf Output:} SOS matrices $S_i$ providing a generalized sum of squares decomposition of $\tilde f$. 
 
In the following example the single steps of the algorithm are illustrated. 
 
\begin{example} 
We consider the symmetric group $S_3$ operating by permutation of the
variables $\{x,y,z\}$.
The group $D_3=S_3$ has order six, with three irreducible representations  
of dimensions $1,1,$ and $2$ (verifying $1^2+1^2+2^2=6$), shown in  
Table~\ref{tab:s3rep}. All three are absolutely irreducible.
 
The invariant ring $\R[x,y,z]^{S_3}$ is the polynomial ring
generated by the elementary symmetric polynomials 
$e_1=x+y+z, e_2=xy+yz+xz,e_3=xyz$, which constitute primary invariants 
(and moreover, a SAGBI basis).

The modules are generated by
\[
b^2_1=(x-y)(x-z)(y-z),\qquad 
b^3_1=\left(\begin{array}{c}
(2x-y-z)/\sqrt 2\\
(y-z)\sqrt 3/\sqrt 2
\end{array}\right),
\qquad b^3_2= 
\left(\begin{array}{c}
( 2yz-zx-xy)/\sqrt 2\\
( zx-xy)\sqrt 3/\sqrt 2
\end{array}\right).
\]
Computing the matrices $\Pi_i$ corresponding to each irreducible  
representation, we obtain:  
\[  
\Pi_1 = 1, \qquad  
\Pi_2 = e_1^2 e_2^2 - 4 e_2^3 - 4 e_1^3 e_3   
+ 18 e_1 e_2 e_3 - 27 e_3^2 , \qquad  
\Pi_3 =   
\left[\begin{array}{cc}  
2 e_1^2 - 6 e_2 & -e_1 e_2 + 9 e_3 \\  
-e_1 e_2 + 9 e_3 & 2 e_2^2 - 6 e_1 e_3  
\end{array}\right].  
\] 
This completes Algorithm I collecting the results which depend on the
group action only.

Now we apply Algorithm II to the global minimization of the specific
quartic trivariate polynomial presented in the Introduction in
Example~\ref{ex:berndpoly}, and previously analyzed in \cite{PaStu01}:
\[  
f(x,y,z) = x^4 + y^4 + z^4 - 4 \, x y z +x + y + z. 
\]

\begin{table} 
\begin{center}  
{\small  
\begin{tabular}{c|c|c|c|c|c|c}  
  & $xyz$ & $yzx$ & $zxy$ & $xzy$ & $zyx$ & $yxz$ \\ \hline  
$\vartheta_1$ &  1 & 1 & 1 & 1 &1 & 1 \\  
$\vartheta_2$ &  1 & 1 & 1 & $-1$ & $-1$ & $-1$ \\  
$\vartheta_3$ &    
$\left[\begin{array}{cc}  
1 & 0 \\ 0 & 1  
\end{array}\right] $   
&   
$\left[\begin{array}{cc}  
-\alpha & -\beta \\ \beta & -\alpha  
\end{array}\right] $   
&   
$\left[\begin{array}{cc}  
-\alpha & \beta \\ -\beta & -\alpha  
\end{array}\right] $   
&   
$\left[\begin{array}{cc}  
1 & 0 \\ 0 & -1  
\end{array}\right] $   
&   
$\left[\begin{array}{cc}  
-\alpha & -\beta \\ -\beta & \alpha  
\end{array}\right] $   
&   
$\left[\begin{array}{cc}  
-\alpha & \beta \\ \beta & \alpha  
\end{array}\right] $   
\end{tabular}  
}  
\end{center}  
\caption{Irreducible orthogonal representations of the symmetric group  
$S_3$. The values of $\alpha, \beta$ are $\alpha = \frac{1}{2}, \beta  
= \frac{\sqrt{3}}{2}.$ Notice that $\alpha^2 + \beta^2 =1,  
\alpha^2-\beta^2 = -\alpha$.}  
\label{tab:s3rep}  
\end{table}  

As shown in \cite{PaStu01}, without using symmetry the associated
semidefinite program has dimensions $10 \times 10$, and the
corresponding subspace has dimension equal to $20$ (codimension is
$35$).

Following Algorithm II the $S_3$-invariant polynomial $f$ is rewritten
as
\[ 
\tilde f({\bf e})= 
e_1^4-4e_1^2e_2+2e_2^2+4e_1e_3-4e_3+e_1. 
\] 
The polynomial $\tilde f$ has weighted degree $4$. Since $\Pi_1$ has
degree zero, the corresponding $S_1$ should have degree four.  We
collect all monomials of weighted degree $2$ in the vector $w_1({\bf
  e}) = (1, e_1, e_2, e_1^2 )$.  The matrix $\Pi_2$ has weighted
degree $6$ and thus does not contribute to the SOS representation of
$\tilde f$.  The entries of $\Pi_3$ have degrees $2,3,4$ resulting in
$w_3(\mathbf{e})=(1,e_1)$.  But the last column and row of $Q_3^{12}$ and
$Q_3^{22}$ are zero since the weighted degree of $\pi^3_{22}$ is four.

Therefore, the original $10\times 10$ semidefinite program is now  
reduced to two smaller coupled SDPs, of dimensions $4 \times 4$ and $3  
\times 3$. The number of variables decreased from twenty to  
five.  

Solving the coupled SDPs numerically, we obtain the value $f^{sos} =
-2.112913882$, that coincides (in floating point arithmetic) with the
algebraic number of degree 6 that is the optimal minimizer $f^*$ (see
\cite{PaStu01} for its minimum polynomial). For easy verification, we
present here a rational approximate solution to the optimal value,
obtained by rounding. This has been done in such a way to preserve the
sum of squares property, thus providing a valid lower bound on the
optimal value.  For this, consider:
\begin{eqnarray*}
S_1 &=& \textstyle\frac{2113}{1000}+e_1+\textstyle\frac{79}{47} e_2 - 
\textstyle\frac{79}{141}e_1^2-\textstyle\frac{1120}{11511} e_1 e_2 -
\textstyle\frac{148}{1279} e_1^3 + 
\textstyle\frac{1439}{2454} e_2^2 -
\textstyle\frac{85469}{188958} e_1^2 e_2+ 
\textstyle\frac{85}{693}e_1^4 \\
&=&
\left[
\begin{array}{c}
1 \\ e_1 \\ e_2 \\ e_1^2 
\end{array}\right]^T
\left[
\begin{array}{cccc}
\frac{2113}{1000} & \frac{1}{2} &  \frac{79}{94} & -\frac{233}{496} \\
\frac{1}{2}&  \frac{13261}{34968} & -\frac{560}{11511} & -\frac{74}{1279} \\
\frac{79}{94} & -\frac{560}{11511} & \frac{1439}{2454} &
-\frac{85469}{377916}\\
-\frac{233}{496} &  -\frac{74}{1279} &  -\frac{85469}{377916} & \frac{85}{693}
\end{array}
\right]
\left[
\begin{array}{c}
1 \\ e_1 \\ e_2 \\ e_1^2 
\end{array}\right], \\
S_3 &=& 
\left[
\begin{array}{cc}
\frac{79}{282}+\frac{74}{1279} e_1 + \frac{304}{693} e_1^2 &
-\frac{2}{9}+\frac{749}{1636}e_1 \\
-\frac{2}{9}+\frac{749}{1636}e_1 &
\frac{3469}{4908}
\end{array}
\right].
\end{eqnarray*}
It is immediate to check that $S_1, S_3$ are indeed SOS, and that they
satisfy $f + \frac{2113}{1000} = S_1 + \langle S_3, \Pi_3 \rangle$,
therefore serving as a valid algebraic certificate for the lower bound
$-2.113$.
\hfill$\diamondsuit$
\label{ex:rational}
\end{example}

In several of our examples, some irreducible representations do not
appear in the final SDP (the corresponding basis is empty). This
happens only because the polynomials are of relatively low degree. It
can be shown, using the Molien series, that \emph{all} the irreducible
representations appear non-trivially in the induced representation on
monomials of total degree $d$, for $d$ sufficiently large (see
\cite{Sta79}).

\section{Computational savings}  
\label{sec:savings}  

The exact amount of computational resources that can be spared at the
SDP solving stage can only be quantified precisely after defining the
specific symmetry group and the way it acts on $\R^n$. In what
follows, we explicitly evaluate the computational savings obtained by
exploiting symmetry, for two specific group families.
  
As mentioned in \cite{PaStu01}, the size of the semidefinite  
programs, when no symmetry is exploited, is equal to  
\[  
N = \binom{n+d}{d}, \qquad m = \binom{n + 2 d}{2d}   
\]  
where $N$ is the dimension of the matrix and $m$ is the codimension  
of the subspace $\mathcal{L}$.  
  
The first example will correspond to polynomials with even 
degrees, invariant under a natural action of $C_2^n$. Secondly, 
we analyze the case of dense polynomials in $n$ variables, of total 
degree $2d$, that are symmetric under the full permutation group 
$S_n$. 
 
\subsection{Even forms} 
  
  In this section we analyze a class of homogeneous polynomials in $n$  
  variables arising from the study of copositive matrices, see  
  \cite{pablo}.  These  polynomials have the symmetry property:  
\begin{equation}  
f(x_1,\ldots,x_n) = f(\pm x_1,\ldots, \pm x_n),   
\label{eq:copocond}  
\end{equation}  
for every possible combination of the signs. The associated group is  
the direct product of $n$ copies of $C_2$, and its order is equal to  
$2^n$.  This group has $2^n$ one-dimensional irreducible  
representations, and a complete set of fundamental invariants is given by  
$\theta_1=x_1^2, \ldots , \theta_n=x_n^2$.  
  
\begin{table}  
\begin{center}  
\begin{tabular}{c|c|c}  
 & $id$ & $s$ \\ \hline  
$\vartheta_1$ (trivial)     &  1 & 1 \\  
$\vartheta_2$ (alternating) &  1 & $-1$  \\  
\end{tabular}  
\end{center}  
\caption{Irreducible unitary representations of the group $C_2$.}  
\label{tab:C2rep}  
\end{table}  
  
Each irreducible representation corresponds to the tensor product of  
$n$ copies of the two representations of $C_2$ (namely, the trivial  
and the alternating), see Table~\ref{tab:C2rep}. Let the \emph{type}  
of a representation of $(C_2)^n$ be the number of factors equal to the  
alternating representation. Then, the Molien series for all $2^n$  
isotypic components can be simultaneously obtained from the identity:  
\[  
\frac{1}{(1-\xi)^n} =   
\sum_{r=0}^n \binom{n}{r}   
\psi_r(\xi), \qquad  
\psi_r(\xi) =    
\frac{\xi^r}{(1-\xi^2)^n} =   
\sum_{d=0}^\infty \binom{n+d-1}{d} \xi^{2d+r},  
\]  
where $\psi_r$ is the series corresponding to each of the $\binom{n}{r}$  
irreducible representations of type $r$. For instance, for the  
invariant ring (i.e., the isotypic component corresponding to the  
trivial irreducible representation) the associated series is:  
\[  
\psi_0(\xi) = \frac{1}{(1-\xi^2)^n} =   
\sum_{d=0}^\infty \binom{n+d-1}{d} \xi^{2d}.     
\]  
The $2^n$ matrices $\Pi_i$ in this case are scalars, and given by 
squarefree products of the $\theta_i$.  
 
Specializing the formulas above to the cubic terms (i.e., only the  
ones in $\xi^3$), we see that:  
\begin{itemize}  
\item There are $n$ components of type 1, each with $n$ monomials of  
  the form $\{x_i^2 x_1,\ldots,x_i^2 x_n \}, i = 1,\ldots,n$.  
\item There are $\binom{n}{3}$ components of type 3, each with one monomial, of  
  the form $\{x_i x_j x_k\}, i < j <k$.  
\end{itemize}   
Therefore, the problem of checking a sum of squares condition for a 
sextic homogeneous polynomial satisfying the symmetry conditions 
(\ref{eq:copocond}) can be \emph{exactly} reduced from an SDP of 
dimension $\binom{n+2}{3}$, to a coupled system of $n$ smaller $n 
\times n$ SDPs with $\binom{n}{3}$ linear constraints 
(cf. \cite[Theorem 5.2]{pablo}). 

The savings can be even more dramatic when higher degree forms are
considered. For instance, the standard approach for a homogeneous
polynomial of degree 8 in $n=10$ variables requires an SDP of
dimension $715\times 715$. The symmetry-aware methodology reduces this
to coupled SDPs of dimensions $55\times 55$ (1), $10 \times 10$ (45)
and $1\times 1$ (210), which are considerably easier to solve. 
The corresponding results for arbitrary order forms can be similarly
obtained from the Molien series above.

\subsection{Symmetric group $S_n$}  
  
Here we analyze the case of polynomials in $n$ variables, invariant
under the symmetric group $S_n$ acting by permutation of the
indeterminates.  As is well known, such $S_n$-invariant polynomials
can be rewritten in terms of the elementary symmetric functions $e_i$.
Since $S_n$ is operating as reflection group the invariant ring is the
polynomial ring $\R[e_1,\ldots,e_n]$. Moreover, the elementary
symmetric polynomials form a SAGBI-basis.
 
We briefly recall some classical results on the representation theory
of the symmetric group $S_n$ \cite{James}. As with all finite groups,
the number of mutually inequivalent irreducible representations is
equal to the number of similarity classes. For $S_n$ there is a
bijection between similarity classes and additive partitions of $n$.
Therefore, the number of irreducible representations is equal to the
number of partitions of $n$ and, by a celebrated result of Hardy and
Ramanujan, asymptotic to $\frac{1}{4 n \sqrt{3}}e^{\pi \sqrt{2n/3}}$.
The dimensions of the representations are given by the number of
standard Young tableaux associated with a given partition. For
instance, for $S_4$ we can see in Table~\ref{tab:Young} that there are
five irreducible representations, of dimensions 1,3,2,3,1,
respectively.
\begin{table}  
\begin{center}  
\begin{tabular}{c|c|c}  
\mbox{Class} & \mbox{Dimension} & \mbox{Young Tableaux}\\ \hline  
$[4]$ & 1 &   
$\begin{array}{cccc} 1 & 2 & 3 & 4 \end{array} $   
\\ \hline  
$[3,1]$ & 3 &   
$\begin{array}{ccc} 1 & 3 & 4 \\ 2 \end{array} \quad   
\begin{array}{ccc} 1 & 2 & 4 \\ 3 \end{array} \quad   
\begin{array}{ccc} 1 & 2 & 3 \\ 4 \end{array}$   
\\ \hline  
$[2,2]$ & 2 &   
$\begin{array}{cc} 1 & 3 \\ 2 & 4 \end{array} \quad   
\begin{array}{cc} 1 & 2 \\ 3 & 4 \end{array}$  
\\ \hline  
$[2,1,1]$ & 3 &   
$\begin{array}{cc} 1 & 2 \\ 3 \\ 4 \end{array} \quad   
\begin{array}{cc} 1 & 3 \\ 2 \\ 4 \end{array} \quad   
\begin{array}{cc} 1 & 4 \\ 2 \\ 3 \end{array}$  
\\ \hline  
$[1,1,1,1]$ & 1 &  
$\begin{array}{c} 1 \\ 2 \\ 3 \\ 4 \end{array}$ \\ \hline  
\end{tabular}  
\end{center}  
\caption{Similarity classes and number of representations for the  
symmetric group $S_4$.}  
\label{tab:Young}  
\end{table}  

\begin{example}
  A symmetric quadratic form in $n$ variables has the representation
  $f(e_1,e_2) = a \, e_1^2 + b \, e_2$. We derive the conditions for
  $f$ to be a sum of squares (or nonnegative, since in the quadratic
  case these notions coincide). By
  Theorem~\ref{thm:matrixsosrep}, and taking into account the degree
  of $f$, if the form is SOS there exists a representation
\[
f(e_1,e_2) = S_1 \cdot 1 + S_2 \cdot (\Pi_2)_{11}
= (c_1 e_1^2) \cdot 1 + c_2 \cdot [(n-1)e_1^2 - 2 n e_2]. 
\]
where $c_1,c_2$ are nonnegative constants.  This gives rise to the SDP
(actually, just a linear program):
\[
a = c_1 + (n-1) \, c_2, \qquad b = -2n c_2, \qquad c_1, c_2 \geq 0. 
\]
Therefore, the necessary and sufficient conditions for $f$ to be a sum
of squares are $(2 n) a + (n-1) b \geq 0, b \leq 0$. 
\end{example}

In general, there seem to be no explicit formulas for the dimensions
of the isotypic components of the polynomial ring (other than the
Molien series), so either asymptotic results or special cases should
be analyzed instead. In the example below, we focus on the symmetric
group in four variables $S_4$.

\begin{example}  
\label{ex:sottile}  
We consider two examples taken from a paper by Sottile \cite{Sottile}
arising from a conjecture of Shapiro and Shapiro. These polynomials
are invariant under the full permutation group in four variables. The
group $S_4$ has five irreducible representations, of degrees equal to
1,3,2,3 and 1.

In Table~\ref{tab:molien} we present the dimensions of the  
corresponding isotypic components when the group acts on the space of  
monomials of degree $d$. The entries were computed using the  
Molien-Hilbert series.   
\begin{table}  
\begin{center}  
\begin{tabular}{c|ccccccccccccccccc}  
  $d=$ & 0 & 1 & 2 & 3 & 4 & 5 & 6 & 7 & 8 & 9 & 10 & 11 & 12 & 13 & 14 & 15 \\ \hline  
$\vartheta_1$ & 1& 1& 2& 3& 5& 6& 9& 11& 15& 18& 23& 27& 34& 39& 47& 54 \\  
$\vartheta_2$ & 0& 1& 2& 4& 6& 10& 14& 20& 26& 35& 44& 56& 68& 84& 100& 120 \\  
$\vartheta_3$ & 0& 0& 1& 1& 3& 4& 7& 9& 14& 17& 24& 29& 38& 45& 57& 66 \\  
$\vartheta_4$ & 0& 0& 0& 1& 2& 4& 6& 10& 14& 20& 26& 35& 44& 56& 68& 84 \\  
$\vartheta_5$ & 0& 0& 0& 0& 0& 0& 1& 1& 2& 3& 5& 6& 9& 11& 15& 18 \\ \hline   
\mbox{Total} & 1& 4& 10& 20& 35& 56& 84& 120& 165& 220& 286& 364& 455& 560& 680& 816  
\end{tabular}  
\end{center}  
\caption{Dimension of the isotypic components of $S_4$ acting on the  
monomials of total degree $d$, not counting multiplicities. The last  
row is the total dimension of the space, equal to $\binom{n+d-1}{d}$.}  
\label{tab:molien}  
\end{table}  
 
Consider for example the case of a quartic form ($2d=4$) appearing in  
\cite[Example 2.2]{Sottile}:  
\[  
\tilde f(\mathbf{e}) = 16e^2_2 - 48 e_1 e_3 + 192 e_4,  
\]  
where the $e_i$ are the elementary symmetric functions in $(s,t,u,v)$.
According to Table~\ref{tab:molien}, the problem can be reduced to a
simple SDP of dimensions $2, 2$ and $1$. The corresponding
symmetry-adapted vector basis is:
\begin{eqnarray*}  
v_1 &=& ( s^2+t^2+u^2+v^2,\, st + su + sv + tu +tv + uv) \\[0.1cm]  
v_{21} &=& ( u^2+s^2 - v^2-t^2,\, u s - v t )\\  
v_{22} &=& ( u^2+v^2 - s^2-t^2,\, u v - s t )\\  
v_{23} &=& ( s^2+v^2 - t^2-u^2,\, s v - t u )\\[0.1cm]  
v_{31} &=& ( u v + s t - s v - t u ) \\  
v_{32} &=&   
\left(\frac{1}{\sqrt{3}}(u v + s t + s v + t u - 2 v t - 2 u s)\right).   
\end{eqnarray*} 
>From the invariant representation viewpoint introduced in
Section~\ref{sec:invar}, the polynomial described has a very simple
form, as it corresponds exactly to a scalar multiple of the $(1,1)$
entry of the matrix $\Pi_3$. In other words, we have:
\[
\tilde f(\mathbf{e}) = 
\left\langle
\left[\begin{array}{cc}
c & 0 \\ 0 & 0 
\end{array}\right],
\Pi_3
\right\rangle,
\]
the value of $c$ depending on the specific normalization chosen for
$\Pi_3$. It is easy then to find the explicit decomposition in terms
of the original variables:
\[  
f(s,t,u,v) = 12 \bigl(u v + s t - s v - t u\bigr)^2 + 12 \left( \frac{1}{\sqrt{3}}(u v + s  
t + s v + t u - 2 v t - 2 u s) \right)^2.  
\]  

We also applied the methods to the degree $2d=20$ symmetric
polynomial in \cite[Theorem 2.3]{Sottile}.  From the table, we can
easily appreciate the benefits of symmetry exploitation: for this,
instead of solving a single SDP of dimension $286\times 286$, five
coupled smaller ones, of dimensions $44,26,24,23,5$ are employed (even
smaller SDPs can be achieved by exploiting the sparsity pattern).

By appropriately choosing the objective function of the SDP, we
quickly found a particular numerical solution with support in only two
of the five isotypic components ($\vartheta_1$ and $\vartheta_2$),
therefore verifying its nonnegativity. Even though Sottile was able to
find an explicit set of polynomials with integer coefficients whose
squares add up to $f$, in the representation we found these do not
seem to have ``nice'' rational expressions. While there is no reason
to expect our methods to yield solutions with particular properties
such as rationality, this is sometimes the case for extremal
solutions. Nevertheless, verification of the solution can still be
done independently, by approximating an interior solution with
rational numbers as was done in Example~\ref{ex:rational}.
\hfill$\diamondsuit$
\end{example}

{\bf Acknowledgments:} The authors would like to thank Bernd Sturmfels
for initiating the cooperation and for helpful discussions, and Kazuo
Murota for the reference to \cite{Murota}.  KG thanks UC Berkeley for
the kind hospitality during her research stay during the spring 2001
and the DFG for support by a Heisenberg scholarship.

\bibliographystyle{plain}  
  

\end{document}